\documentclass[11pt]{amsart}

\usepackage{amssymb} 
\usepackage[mathscr]{eucal} 
\usepackage[all]{xy}
\usepackage{amsmath} 

\def\cal#1{{ \mathcal {#1} }} 
\def\frak#1{{ \mathfrak {#1} }} 

\newcommand{\tagnum}[1] 
{\noindent \hbox to 22pt {\hss\bf #1.}}

\theoremstyle{plain}

\newtheorem{thm}{Theorem}[section]
\newtheorem{cor}[thm]{Corollary}
\newtheorem{lem}[thm]{Lemma}
\newtheorem{prop}[thm]{Proposition}
\newtheorem{defprop}[thm]{Definition-Proposition}

\newtheorem{prob}[thm]{Problem}

\newtheorem{defn}[thm]{Definition}

\newtheorem{rem}[thm]{Remark}

\numberwithin{equation}{section}

\def\lar{\leftarrow}
\def\rar{\rightarrow}
\def\rxar#1{\xrightarrow{\kern -1pt #1}}
\def\lxar#1{\xleftarrow{#1 \kern -3pt}}

\newcommand{\CZ}{{\mathbb{C}}}
\newcommand{\FZ}{{\mathbb{F}}}

\newcommand{\ZZ}{{\mathbb{Z}}}

\newcommand{\Pn}[1]{{ \mathbb{P}^{#1}}\relax}

\newcommand{\fr}[2]{{#1}/{#2}\, }

\newcommand{\QED}{\hspace*{\fill}$Q.E.D.$}

\newcommand{\QQ}{\ensuremath{\mathbb{Q}}}

\newcommand{\ra}{\ensuremath{\rightarrow}}

\def\eea{\end{eqnarray*}}
\def\bea{\begin{eqnarray*}}

\def\A{{\mathcal{A}}}

\def\C{{\mathbb{C}}}
\def\E{{\mathcal{E}}}
\def\F{{\mathcal{F}}}

\def\L{{\mathcal{L}}}

\def\P{{\mathcal{P}}}

\def\Z{{\mathbb{Z}}}
\def\hol{{\mathcal{O}}}
\def\PP{{\mathbb{P}}}
\def\i{{\iota}}

\DeclareMathOperator{\rk}{rk}
\DeclareMathOperator{\cork}{corank}

\newcommand{\shExt}{\cal Ext\kern 2pt}
\newcommand{\shHom}{\cal Hom\kern 2pt}
\newcommand{\Syz}{\cal Syz\kern 2pt}

\DeclareMathOperator{\coker}{coker}

\newcommand{\Proof}{{\it Proof. }}

\def\hd#1{\hbox to 19pt{\hfil$#1$\hfil}} 

\def\hdotrulefill{\cleaders\hbox
{\hbox to 2pt{\cleaders\hrule\hfill}\kern 2pt}
\hfill}

\def\vdotrule{\hbox to 0 pt{\hfill
\vbox to 10.5pt{\kern 1pt
\hbox to 0pt{\hfil \vrule height2pt \hfil} \kern 2pt
\hbox to 0pt{\hfil \vrule height2pt \hfil} \kern 2pt
\hbox to 0pt{\hfil \vrule height2pt \hfil} \kern 2pt
\hbox to 0pt{\hfil \vrule height2pt \hfil} \kern 2pt
\vss}\hfill}}

\title{Even sets of nodes on sextic surfaces}

\author{Fabrizio Catanese}

\address{ Lehrstuhl Mathematik VIII,
Mathematisches Institut der
Universit{\"a}t  Bayreuth, NW II \\
D-95440 Bayreuth,  Deutschland.}

\author{Fabio Tonoli }

\email{fabrizio.catanese@uni-bayreuth.de}
\email{fabio.tonoli@uni-bayreuth.de}

\date{\today}

\begin{document}

\begin{abstract}
We determine the possible even sets of nodes on sextic surfaces in $\Pn 3$,
showing in particular that their cardinalities are exactly the numbers
in the set $\{ 24, 32, 40, 56 \}$.
We also show that all the possible cases admit an explicit description.
The methods that we use are an interplay  of coding theory and projective
geometry on one hand, of homological and computer algebra on the other.

We give a  detailed geometric construction for the new case
of an even set of 56 nodes , but the ultimate verification of existence
  relies on computer calculations.
Moreover, computer calculations
have been used  more than once
in our research in order to get good guesses.

The construction gives a maximal family, unirational and
of  dimension 27,
of nodal sextics with an even set of 56 nodes.

As in \cite{CaCa} (where the other cases were  described)
  each such nodal surface $F$  is given as the determinant of a symmetric
map $ \varphi : \E^\vee \ra \E$, for  an appropriate vector bundle $\E$
depending on $F$.
The first  difficulty  here is to show the existence of
such vector bundles. This leads us to the
investigation  of a hitherto unknown
  moduli space  of rank 6 vector bundles which we
show to be  birational to a moduli space
of plane representations of cubic surfaces in $\PP^3$.
The resulting picture shows a very rich and interesting geometry.
The main difficulty is to show the existence of "good" maps $\varphi$,
and the interesting phenomenon which shows up is the following:
the "moduli space" of such pairs $(\E, \varphi )$ is
(against our initial hope) reducible,
and for a general choice of $\E$ the determinant of
  $ \varphi $ is the double of a cubic surface $G$.
Only when the vector bundle $\E$
corresponds to a reducible cubic surface, then
we get an extra component of the space
of such pairs $(\E, \varphi )$, and a general choice
in this component yields one of our desired
nodal sextic surfaces.

\end{abstract}

\maketitle

\section*{Introduction}

Let $F$  be a nodal surface in $\Pn 3$ of degree $d$: i.e., $F$  has only
$\mu$
nodes (ordinary double points)
$P_1, \dots P_{\mu}$ as singularities.

A natural and classical question is
to ask for the maximum possible number of nodes $\mu(d)$ that such a
surface $F$ can
have.

The theory of projectively dual surfaces shows easily that $ \mu(d) <
\frac{1}{2} d (d-1)^2$ for $d \geq 3$
and the slightly better inequality given by Bassett in 1907
(cf. \cite{bass}) was obtained using this method.

The function $\mu(d)$ is only known for $ d \leq 6$,
and for $ d \leq 5$ one has an explicit description of
the nodal  surfaces which attain the maximum $ \mu(d)$:
the Cayley cubic, the Kummer quartics, and the Togliatti
quintics (cf.  \cite{cay1}, \cite{cay2}, \cite{kum},
\cite{tog1}, \cite{tog2}, \cite{Bea}).

An important tool to investigate the function $\mu(d)$
for small values of $d$ ($ d \leq 17$), and to characterize the
maximizing surfaces was  introduced by
Beauville in \cite{Bea}: he attached
a binary code to each nodal surface $F$
and used coding theory in order to show that $\mu(5)=31$.

The method of using coding theory was later used by Jaffe and Ruberman
in order to show (see \cite{JR}) that $\mu(6) = 65$,
but their proof is not so short as the one by Beauville,
partly because at that time a complete knowledge
about the cardinality of an even set of nodes\footnote{
 We adopt  here the terminology of \cite{CaCa} concerning
the notion of even sets of nodes which was introduced in \cite{Ca}:
namely, the strictly even sets of \cite{Ca} are called even sets, while
the weakly  even sets of \cite{Ca} are called half-even, or 
$1/2$-even set of nodes.}  on a sextic was missing (the binary code consists
of  the even sets of nodes on $F$,
introduced in \cite{Ca}, where a complete classification of even sets for
  degree $ d = 5$ was given).

Today we still ignore if the Barth 6-ics (see \cite{Ba}) are those 
which achieve
the maximum $\mu(6)= 65$
and until now an explicit description
of the possible even sets of nodes for sextic surfaces was missing.
A general structure theorem for  even and $1/2$-even 
sets
was given in \cite{CaCa}: 
but the cases where the cardinality $t$ of an even set
would be  $ > 40$ were  excluded only as a consequence of
a conceptual error which was pointed out to the authors by 
Duco van Straten.
Thus the simple proof by J. Wahl (cf. \cite{wahl}) of $\mu(6)= 65$ also became
invalid.

We rescue the situation here by showing the following

{\bf Main Theorem A.}
{\em  Let $F$ be a  nodal surface of degree $d=6$ in $\Pn 3$ with an
even set of $t$ nodes. Then $t\in \{24,32,40,56\}$.
These four possibilities  occur and can be explicitly described.}

The situation is thus more complicated than for $d \leq 5$,
  the list of possible
cardinalities $t$ is (cf. e.g. \cite{CaCa}):
$$
\begin{tabular}{cl}

$d=3$ & $t=4$ \\
$d=4$ & $t\in\{8,16\}$ \\
$d=5$ & $t\in\{16,20\}$ \\
$d=6$ & $t\in\{24, 32, 40, 56 \}$ \\
\end{tabular}
$$

We  first show  in section 1 that the case of an even set of $64$ nodes cannot
exist. The simple new idea is to study  the so called extended code (cf. e.g.
\cite{Ca2}) and we then use a mixture of geometric
and coding theory arguments, as was done in the papers cited above,
for instance in \cite{JR}, where  the
case of an even set of $48$ nodes  was excluded .

We  then proceed, using the structure theorem of \cite{CaCa},
to construct  explicit cases of sextics with an even set of 56 nodes.

The bulk of the paper is devoted to this purpose,
and we get the following result.

{\bf Main Theorem B.}
{\em There is a  family of nodal sextic surfaces
  with $56$ nodes forming an even set, parametrized by  a
smooth irreducible rational variety $\Phi_0$ of dimension
33, whose image $\Xi_0$ is a unirational subvariety  of dimension 27
of the space of sextic surfaces.
Moreover, the above family is versal, thus $\Xi_0$ yields
an irreducible  component of the subvariety of  nodal sextic surfaces with 56
nodes.}

\medskip

The fact that  a maximal family
of nodal sextics with 56 nodes forming an even set 
has dimension equal to  27
means that these nodes impose
  independent conditions on the space of sextic surfaces
(cf. \cite{b-w}).
It is an interesting question to find the smallest degree for
which there exist even sets of nodes failing to impose independent
conditions. 

\medskip

As already mentioned, it follows from the more general result
of  \cite{CaCa} that
every even  nodal set on a sextic surface $F$ occurs as the corank  2
degeneracy  locus of a symmetric
map $ \varphi : \E^\vee \ra \E$, for  an appropriate vector bundle $\E$
depending on $F$.

The method to construct $\E$ is based on a combination
of Beilinson's theorem and of a revisitation of
Horrocks' correspondence due to Charles Walter (cf.
\cite{wal}),
which was exploited in \cite{CaCa}.
The bundle $\E$ is constructed starting from a
submodule
$M$ of the intermediate cohomology module $H^1_* (\F)$
of the quadratic sheaf $\F$ associated to the even set, 
and corresponding to the choice of a Lagrangian
subspace $U$ of $H^1 (\F(1))$.
The choice of $M$ determines a unique vector bundle $\E$, if
a certain generality assumption (which we call {\em first assumption})
is verified.

The construction is quite explicit if we make another generality
assumption, namely that the two nonzero degree components
of the artinian graded module $M$,
the previously mentioned $U$ and another one denoted by $W$,
both have dimension equal to 3.
If we denote as customary by $V$ the vector space of linear forms on $\PP^3$,
then the module $M$ is completely determined
by the multiplication tensor $B \in  U^\vee \otimes  V^\vee \otimes W$
for $M$.

We then show that the tensor $B$ determines explicitly
the bundle $\E$ as the kernel of an exact sequence
$$ 0 \ra \E \ra U \otimes V \otimes \hol(1) \ra (W \otimes \hol(1))
\oplus ( U \otimes \hol(2)) \ra 0$$
where the first component is precisely $B$, and the second is the
standard Euler map, here denoted by  $\epsilon$.

Section 3 then ends by showing that the family of pairs $(\E,\varphi)$
is parametrized (non uniquely) by the following family of pairs
\begin{equation*}
\begin{split}
\frak M_{AB} : = &\{(B,A)\mid
B\in  U^\vee \otimes  V^\vee \otimes W,\\&
A\in ( U\otimes V)\otimes ( U\otimes V) \otimes H^0(\cal O_\Pn 3(2)),
A=^t\!A, (B,\epsilon) \cdot A=0\}.
\end{split}
\end{equation*}
$\frak M_{AB} $ sits inside an affine space of dimension $ 816$,
and it is not possible even for the computer to find the
decomposition of $\frak M_{AB} $ into irreducible components.
It is clear that $\frak M_{AB} $  dominates the space of the above
tensors $B$, and, if $\frak M_{AB} $ were irreducible, one would
obtain the sextic surfaces immediately by a random choice.

However, for long time all the random choices would always
give the double of a cubic surface $G$ as determinant of $\varphi$,
and it looked like even sets with 56 nodes would not exist.
We then tried to prove that this was indeed the case,
and we had to find an explanation for
the cubic surface $G$.

Now, it is classical that to a $ 3 \times 3 \times 4$ tensor $B$ one 
can associate
a cubic surface in $\PP^3$ by taking the determinant
of the corresponding  $ 3 \times 3 $ matrix of linear forms on $\PP^3$.
However, in our case we get a cubic surface $G_*$  in the dual projective space
${\PP^3}^\vee = Proj (V^\vee)$, together with two
different realizations of $G^*$  as a blow up of a projective
plane $ Proj (U^\vee) $ (respectively, $ Proj (W) $) in a
set of 6 points. These are the points where the
Hilbert-Burch $ 3 \times 4$
matrix of linear forms on $ U$ drops rank by 1,
and the rational map to ${\PP^3}^\vee$ is given
by the system of cubics through the 6 points,
system which is generated by the determinants of the 4  $3 \times 3$
minors of the Hilbert Burch matrix.

One passes from one realization to the other simply by transposing
the tensor, and we will call this the trivial involution
for  $ 3 \times 3 \times 4$ tensors: but what we have discovered,
through geometry,
is the existence of another involution for
$ 3 \times 3 \times 4$ tensors, which
we called the {\bf cross-product involution}.

This second involution associates to a general tensor
$B\in  U^\vee \otimes  V^\vee \otimes W$ another tensor
$\cal B\in {W'}^\vee\otimes V\otimes U'$,
where $ W' := \Lambda^2 (W)$ and $U'$ is defined as the kernel
of the map $\Lambda^2(W^\vee)\otimes V \ra  U^\vee
\otimes W^\vee $
  induced by contraction with $B$ (cf. \ref{cross} for the proof that
we have indeed a birational involution).

In fact, to $\cal B$ corresponds now a cubic surface $G \subset \PP^3$,
which is related to a general bundle $\E$ through the existence
of an exact sequence
$$  0\rar 6\cal O \rar \cal E \rar \tau \rar 0,$$
where $\tau$ is an invertible sheaf on the cubic surface $G$.
One can see more precisely that $\cal B$ determines
a sheaf $\cal G$ on $G$ such that $\tau = \cal G^{\otimes2} (-1).$

We found in this way a nice explanation of the phenomenon pointed out by
the computer: we got as determinant the surface $G$ counted twice,
simply because,  in view of the above exact sequence,   for a smooth cubic
surface (indeed, irreducible) all the symmetric
endomorphisms
$ \phi \in H^0 (S^2 (\E))$ are induced by the inclusion
$  S^2 (6 \hol) \ra  S^2
(\E)$.

It was clear at this point that if $ H^0 (S^2 (\E))$ would always 
have dimension
21, then we could not get any nodal sextic surface of the desired type,
but it was of course possible that the dimension could jump up
for special surfaces $G$, and that our parameter space $\frak M_{A,B}$
would be reducible. As explained in section 7, a small computational
simplification and the reduction to finite fields allowed
to make many more random attempts, until the first sextic surface
appeared. Since a determinantal approach  predicts
that the space of tensors $B$ for which the dimension
of  $ H^0 (S^2 (\E))$ jumps has codimension 7, it was only natural
to guess that the case which works out is the case of
tensors $\cal B$ corresponding to  reducible cubic surfaces.
This guess turned out to be true.

The cross product involution can also be phrased as a duality
theorem for a certain moduli spaces of vector bundles on $\PP^3$.
Namely, we prove the following

{\bf Theorem C}
  {\em    Consider the moduli space $ \frak M^s ( 6; 3,6,4)$ of  simple
rank 6  vector bundles $\cal E$ on $\Pn 3$
     with Chern polynomial $ 1 + 3t + 6 t^2 + 4 t^3$,
and inside it the  open set $\frak A$  corresponding to the simple 
bundles with minimal cohomology, i.e., such that
\begin{enumerate}
\item
$ H^i(\E) = 0 \ \forall i \geq 1$
\item
$ H^i(\E (-1)) = 0 \ \forall i \neq 1$
\item
$ H^i(\E (-2)) = 0 \ \forall i \neq 1$.
\end{enumerate}

Then   $\frak A$  is irreducible of dimension 19
and it is bimeromorphic to $ \frak A^0 $,
where $\frak A^0$ is an open set of the G.I.T. quotient space of the
projective space $\frak B$ of
tensors of type $(3,4,3)$, $\frak B : = \{B\in \PP ({U}^\vee\otimes
V^\vee\otimes W)\}$
by the natural action of $SL(W) \times SL(U)$.

Let  moreover $ [B] \in \frak A ^0$ be a general point:
then to $[B]$
corresponds a vector bundle $\E_B$ on $\PP^3$, and also a vector bundle
$\E^*_B$ on ${\PP^3}^\vee$, obtained from the direct construction
applied to $\cal G^*_B$ (cf. \ref{*}).  $\E^*_B$
is the vector bundle $\E_{\cal B}$, where $[{\cal B}] \in {\frak A ^0_*}$
is obtained from $B$ via the cross product involution.}

\bigskip

It would be interesting to further investigate the moduli space
of  Gieseker semistable  rank 6 vector bundles
$\cal E$ on $\Pn 3$
     with Chern polynomial $ 1 + 3t + 6 t^2 + 4 t^3$.
  In section 5 we make a first step in this direction
proving the Mumford-Takemoto
semistability  of the general bundle in $\frak A$.

Section 7 is devoted to a brief account of the random approach
which we already mentioned, while the appendix contains the two
Macaulay scripts which are needed for the ultimate verification of
the existence
of surfaces which have an even set of 56 distinct  nodes as 
the only singularities.


\section{Excluding via coding theory}
\label{CodingTheory}

Throughout this section $F$ will be a normal surface in $\Pn 3$ of
degree $d$ having
at most Rational Double Points as singularities, and possessing moreover $\mu$
nodes (ordinary double points)
$P_1, \dots P_{\mu}$ among its singularities.

We let $\pi:\tilde F\rar F$ be the minimal resolution of the
singularities of $F$.
It is well known (cf. \cite{tju}) that $\tilde F$ is diffeomorphic to
a smooth surface of degree
$d$ in $\Pn 3$: in particular $\tilde F$ is simply connected and for its
second Betti number we have $ b_2(\tilde F) = d (d^2 - 4 d + 6) -2$.

We let  $A_1,\ldots,A_\mu$ be the exceptional $(-2)$-curves ($\cong \Pn
1$) coming from the
blow up of the nodes $P_1, \dots P_{\mu}$, and we let $H$ be the full
transform of
a  plane section of $F$.

Let $V$ be the $\ZZ /2$-vector space freely generated by the
$A_i$'s and consider the map
$$\epsilon : V :=\oplus_{i=1}^\mu ( \ZZ   /2 ) A_i \rar H^2(\tilde
F,\fr{\ZZ}2),$$ given by the reduction modulo two of the integral first Chern
class of a divisor:
$\epsilon  (\Sigma_i a_i A_i) := c_1( \Sigma_i a_i A_i)
\mod(2)$.
Let $U$ be the image of $\epsilon$.

Since  $A_i\cdot A_j=-2
\delta_{ij}, A_i \cdot H=0, H^2=d$, it follows that $U$ is an isotropic
subspace of $V$, and since the intersection product modulo $2$ is non
degenerate
its dimension does not exceed
$b_2(\tilde F)/2$.

In the case where the surface has even degree $d\equiv 0 (\text{mod } 2)$,
we consider more generally $\tilde V:= V\oplus \fr{\ZZ}2 H$,
        $\tilde {\epsilon} : \tilde V \rar  H^2(\tilde F,\fr{\ZZ}2),$ and the
corresponding isotropic subspace $\tilde U :=  Im ( \tilde {\epsilon})$.

\begin{defn}

        1) The {\bf strict code} $K$ associated to
the nodal set $\{ P_1, \dots P_{\mu} \}$  on the surface $F$ is
the binary code  $ K : = {\rm ker} (\epsilon).$

2) If $d\equiv 0 (\text{mod } 2)$ the {\bf enlarged code} $\tilde K$
associated to
the nodal set $\{ P_1, \dots P_{\mu} \}$  on the surface $F$ is
the binary code  $ \tilde  K : = {\rm ker} (\tilde  \epsilon).$

\end{defn}
By the above inequality for  $ dim (U)$, we get
$$\dim K \geq \mu   -
\frac{1}{2} b_2(\tilde F) = \mu - \frac{1}{2} d (d^2 - 4 d + 6) + 1,
\dim \tilde K \geq  \mu  -
\frac{1}{2} d (d^2 - 4 d + 6) + 2.$$

\begin{rem}
By Miyaoka's inequality (cf. \cite{miya}) $ \mu \leq \frac{4}{9} d (d-1)^2$,
therefore only for $ d \leq 17$ we get for sure a nontrivial code $K$ :
because $ {\rm dim } K \geq
\mu - \frac{1}{2} d (d^2 - 4 d + 6) + 1.$
\end{rem}

Notice that the notion of an even, respectively half-even, set of nodes
can be derived from the coding-theory framework.

\begin{defn}\label{notation}

A vector $v\in V$ is completely determined by  its
{\bf support} $N_v:=\{i \mid v_i=1\}$.
The cardinality of the support is called the {\bf weight} of $v$
and denoted by $w(v):=\#N_v$.

By the universal coefficients theorem and Lefschetz' $(1,1)$ theorem
the condition $v\in
K$ is equivalent to the 2-divisibility of $\sum_{i\in N_v} A_i$ in
$Pic(\tilde F)$.
We denote  by $L$ a divisor on $\tilde F$ such that $ 2 L \equiv
\sum_{i\in N_v} A_i$.
The class of $L$ in $Pic(\tilde F)$ is uniquely determined, because
$Pic( \tilde F)$ has no torsion.
We have then a finite double cover
$\tilde S$ of $\tilde F$ branched exactly on the nodal curves $A_i$
such that $i\in N_v$,
and moreover $ f_* \hol_{\tilde S} = \hol_{\tilde F} \oplus \hol_{\tilde F}
(-L).$ Correspondingly, we have a double cover $ f : S \ra F$, with
$ f_* \hol_S = \hol_F \oplus \F$, ramified
exactly in $\Delta:=\{P_i\in F \mid i\in N_v\}$ (cf. \cite{Ca} and
\cite{CaCa} for more
details).
\\
These sets $\Delta$ are called {\bf even sets of nodes} (cf. \cite{Ca}).

\medskip
Similarly, one defines a {\bf  half-even set of nodes} $\Delta$
by the condition that its associated word $\tilde
v:=(v_1,\ldots,v_\mu,1)$,
obtained by setting $v_i=1 \iff P_i \in \Delta$,
belongs to the enlarged code $\tilde K$. 
This condition is again
equivalent to the
existence of a divisor $L$  in $Pic(\tilde F)$ with $ 2 L \equiv
\sum_{i\in N_v} A_i +H$.

We define the {\bf weight} and {\bf support} of $\tilde v$
as the weight and support of the word $v:=(v_1,\ldots,v_\mu)\in V$
(these notions are different from the corresponding ones in coding theory).

\medskip
Observe finally that $K = \tilde K \cap \{ \tilde v | \tilde v_{\mu +
1} = 0 \}$.
\end{defn}

As shown in  \cite{Ca}, Prop. 2.11 and Prop. 2.13,
the geometric interpretation of even sets of nodes
in terms of  double coverings
allows to give the following restrictions for
the cardinality $t$ of an even (resp.: half-even) set of nodes

\begin{prop} \label{2.13Babbage} Let $t : = w(u)$ be the weight of
a code word $u$.
         \begin{enumerate}
         \item Assume $u\in K$ : then $t\equiv 0 \ (4)$.
           Moreover, if $d$ is even, then $t\equiv 0 \ (8)$.
         \item If ($d$ is even) $u\in \tilde K\setminus K$, then   $t\equiv
d(2d-7)/2 \ (4)$. In particular, for $d=6$, $t\equiv -1 \ (4)$.
         \end{enumerate}
\end{prop}

\begin{cor}\label{code}
Let $d = 2 (2k+1)$ be twice an odd integer and assume that $K, \tilde K$
are the codes corresponding to an even set of nodes $\Delta$:
then $\tilde K = K$.
\end{cor}
\Proof
Our  assumption is that the code $K \subset \tilde K$ contains the vector
$\mathbb I$ whose
coordinates are all equal to $1$, except the last which equals $0$.
If we have a vector $ w \in \tilde K \setminus K$ and let $t$
      be its weight, then
the weight of $ ( \mathbb I + w) \in \tilde K \setminus K$
is congruent to $- t$ modulo $(4)$. Since  $t\equiv d(2d-7)/2
\equiv 2 k +1 \ (4)$,
then $-t \not \equiv d(2d-7)/2 \ (4)$,  contradicting (2) of the
previous proposition.

\QED

Let us examine by means of coding theory  which even
sets of nodes can
occur on sextic nodal surfaces.
The main result of this section is the following theorem.

\begin{thm}\label{noevensetwith64}
On a sextic normal surface $F$ with only Rational Double Points as
singularities
there does not exist an even
set of nodes of cardinality $t=64$.
\end{thm}

In order to prove the theorem we first prove some preliminary results.

\begin{lem}
Suppose that there exists an even  set $\Delta$ of nodes of
cardinality $t = 64$ on
a normal sextic  surface $F$.

Let $\gamma:\ F \dashrightarrow F^\lor \subset \Pn 3$ be the Gauss map of $F$,
        given by the partial derivatives $(\frac{\partial F}{\partial x_i})$.

1) $\gamma$ corresponds to a linear subsystem  $\cal L$ of $| 5 H -
\sum_{i=1}^{t} A_i|$
on $\tilde F$ whose fixed part $\Phi$  is contained in the preimage
of the singular points
of
$F$ which are not the nodes of $\Delta$.

2) Let $L$ be a divisor on $\tilde F$ with $ 2 L \equiv \sum_1^{t} A_i $: then
$ H^0 ( \tilde F, \hol_{\tilde F} ( 2 H - L)) = 0.$
\end{lem}

\Proof
The first assertion follows since the zero locus of the partial derivatives
$(\frac{\partial F}{\partial x_i})$ on $F$ is exactly the singular
locus of $F$, and at
each node the
        partials $(\frac{\partial F}{\partial x_i})$ define the maximal
ideal. Thus $\gamma$
is a morphism around each $(-2)$-curve $A_i$, which is indeed embedded as a
plane conic.

Assertion 2) is proven by contradiction. Assume in fact that  $ C \in
| 2 H - L|$.

We calculate now $ C \cdot (5 H - \sum_{i=1}^{t} A_i ) = 60 - 64 = -4$.

However, $ \Phi \cdot H =  \Phi \cdot A_i = 0$ by our first assertion,
whence $ C \cdot (5 H - \sum_{i=1}^{t} A_i )$ equals the intersection number
of $C$ with the movable part of the linear system
        $\cal L$, which is obviously a non negative number.

We have obtained the desired contradiction.

\QED

\begin{prop}
Suppose there exists an even set of nodes of cardinality $64$ on
a sextic normal surface with only Rational Double Points as singularities,
       and let $ f : S \ra F$ be the
corresponding finite double cover.

       Then $h^1(S,\hol_S)=5$.
\end{prop}

\begin{proof}
We have
$$h^1(S,\hol_S) = h^1(F,\cal F)=h^1( \tilde F, \hol_{ \tilde F} (-L)).$$
Moreover, since $ H \cdot L = 0$ one easily sees that
$h^0( \tilde F, \hol_{ \tilde F} (-L))=0$ and argues then that
$h^2( \tilde F, \hol_{ \tilde F} (-L))=h^0( \tilde F, \hol_{ \tilde
F} (2 H + L))=
h^0( \tilde F,\hol_{ \tilde F} (2H-L)=0$ by the previous lemma,
the second equality following from the fact that every divisor in $|
2H+L | $  contains $\sum_{i=1}^{t} A_i$.

Whence $- h^1(S,\hol_S) = - h^1( \tilde F, \hol_{ \tilde F} (-L))=
\chi( \tilde F, \hol_{ \tilde F} (-L)) =\chi(\hol_{ \tilde F})+
1/2 (-L)\cdot(-L-2H) =11-16=-5$.

\end{proof}

\begin{lem}
Suppose there exists an even set of nodes $\Delta$ of cardinality $64$ on
a sextic normal surface $F$ with only Rational Double Points as singularities.
       Let $K$, $\tilde K$ denote the
corresponding binary  codes.

Then $\dim \tilde K =12> \dim K =11$.
\end{lem}

\begin{proof}

By the previous proposition the surface $S$, the finite double cover
of $F$ ramified
exactly along $\Delta$, has invariants $p_g=10, q=5, K_S^2=48, \chi(S)=6$.

The corresponding (non minimal) smooth surface $\tilde S$, the double cover
of $\tilde F$, has the same invariants as $S$.

By \cite[Lemma 2]{Bea} or \cite[Thm 4.5]{JR} it follows that the code
$K$ has dimension $b_1(\tilde S)+1=11$.
We already remarked that  $ {\rm dim} \tilde K \geq 65 - 53 = 12,$
and it is obvious that $K \subset \tilde K$ has codimension at most $1$.

\end{proof}

\begin{proof}[Proof of Thm. \ref{noevensetwith64}]

The previous lemma contradicts then corollary \ref{code}.

\end{proof}

As an immediate consequence of Thm. \ref{noevensetwith64} and
\cite[Sec. 7]{JR},
we obtain the following.

\begin{cor}\label{weights}
Let $F$ be a a sextic normal surface  in $\Pn 3$ with only Rational
Double Points as singularities
       with an even set
of $t$ nodes. Then $t\in \{24,32,40,56\}$.
\end{cor}

\Proof

Since $ t \equiv 0  \ (8)$, the inequality $ t \leq 64$ follows
      from the classical inequalities of Bassett
and of Miyaoka, and the case $t=64$ was just excluded.
In \cite{CaCa} it is shown that $ t \geq 24$, and that the cases $t =
24, 32, 40$
do exist.

The non-existence of even sets of 48 nodes on nodal sextics is proven
in \cite[Sec. 7]{JR}.

\QED

\section{Cohomology modules and bundle symmetric maps}

In this section, after recalling the main result of \cite{CaCa},
namely the correspondence between even sets of nodes and bundle
symmetric maps, we shall give bounds for the cohomology groups
$H^i(\F(j))$ of the quadratic sheaf $\F$ associated to an even set
of nodes $\Delta$.

\bigskip
We first recall the main result of \cite{CaCa}, according to the
following notation:
$\delta \in \{ 0,1\}$ and $\fr\delta2$-even stands for even if $\delta=0$,
respectively for half even if $\delta=1$.

\begin{thm}[{\cite[Thm. (0.3)]{CaCa}}]\label{CaCaThm}
Let  $\Delta$ be a $\fr\delta2$-even set of nodes on a normal surface
$F$ of degree $d$, let
       $f: S \rar F$ denote a corresponding double cover of
$F$, let $\cal F$ be the
anti-invariant part of $f_* \cal O_S$.

Then there exists a locally free sheaf $\cal E$ on $\Pn 3$ and a symmetric map
$\varphi$ yielding an exact sequence
$$
(**) \ 0 \rar \mathcal E^{\vee}(-d-\delta) \rxar\varphi \mathcal E \rar
\mathcal F \rar 0.
$$

In particular, $F = \{ x\, |\, \det (\varphi (x)) =0 \}$, 
$\Delta =\{ x\, |\,  \cork
(\varphi (x)) \geq 2\}$.

Conversely, assume that one is given an exact sequence as in $(**)$ with
$\varphi$ symmetric, such that  
$F = \{ x\, |\, \det (\varphi (x)) =0 \}$
is a normal surface
and $\Delta : =\{ x\, |\, \cork(\varphi (x)) \geq 2\}$ is a reduced set
of $t$ points:
then $\Delta$ is a $\fr\delta2$-even set of nodes on $F$.
\end{thm}
\noindent

The ideal of the subscheme $\Delta$ is  the second Fitting
ideal of $\varphi$, i.e., on local trivializing affine sets for $\cal E$,
it is given by the determinants of the  $(\rk \cal E -1)$-minors
of $\varphi$.

We briefly explain how
the sheaf $\cal E$ is explicitly constructed in \cite{CaCa} by means of
C. Walter's interpretation of Horrocks'
correspondence (cf. \cite{wal}).

Assume that  the intermediate cohomology module
$H^1_*(F,\cal F):=\oplus_{i\in \ZZ} H^1(F,\cal F(i))$ is known
(it is an Artinian  graded module over the polynomial ring of $\Pn 3$,
$\A : = \CZ [x_0,x_1,x_2,x_3]$).

One considers then the (Artinian) graded module
\begin{equation}\label{M}
M:= U \bigoplus\,  (\oplus_{i>\fr{(d-4)}2} H^1(F, \cal F(i)) ),
\end{equation}
where, if $d$ is even, $U$ is a Lagrangian subspace in the
Serre self-dual cohomology space
$H^1(F, \cal F(\fr{(d-4)}2))$, and $U:=0$ if $d$ is odd.

Recall that the first syzygy bundle $Syz^1(M)$ is obtained
from a projective graded resolution of the module $M$ by free $\A$-modules
$$ 0 \rar \P^4 \rar \dots \rar P^1 \rxar{\alpha_1} P^0
\rxar{\alpha_0} M \ra 0$$
as follows:
the homomorphism $\alpha_1: P^1 \rar P^0$ induces a corresponding
homomorphism $(\alpha_1)^\sim$ between the (Serre-) associated sheaves
$({P^1})^\sim$ and
$({P^0})^\sim$ and the first syzygy
bundle of $M$ is defined as $Syz^1(M) : = {\rm Ker} ( \alpha_1^\sim )$.

One has a natural homomorphism $Syz^1(M) \ra \F$
(cf. \cite{CaCa}, pages 240-1) induced by truncation, whence
one gets a homomorphism $H^0_*(Syz^1(M)) \ra H^0_*(\F)$,
which needs not be surjective.

The bundle $\cal E$ is then defined as  the direct sum of $Syz^1(M)$
with a direct sum of line bundles, whose generators 
induce a minimal set of generators of the cokernel of 
$H^0_*(Syz^1(M)) \ra H^0_*(\F)$,
in order that one obtains
a surjection between $H^0_*(\Pn 3,\cal E)$ and $H^0_*(F,\cal F)$.

Thus a first important step is the one of determining the
intermediate (Artinian) cohomology module
$H^1_*(F,\cal F):=\oplus_{i\in \ZZ} H^1(F,\cal F(i))$,
in particular one has to determine the possible dimensions
of its graded pieces, i.e., the numbers $h^1(\cal
F(i))$. Later on, when we want to impose the surjectivity
of $H^0_*(\Pn 3,\cal E) \ra H^0_*(F,\cal F)$ it will also
be important to determine the dimensions $h^0(\cal
F(i))$.

In short, the first necessary task is to determine the possible values
for the cohomology table $h^j(\cal
F(i))$ of $\cal F$
(a priori just $\chi(\cal F(i))$ is known, and it is determined by
the degree $d$ of $F$ and the number $t$ of nodes of $\Delta$).

Besides geometrical estimates, an important tool used in \cite{CaCa}
is the Beilinson complex (cf. \cite{Bei}) constructed
from the cohomology table $h^j(\cal  F(i))$.

\begin{rem}\label{bei}
It is well known that for any coherent sheaf $\cal G$ on $\mathbb P^n$
the complex
$\cal K^i:=\bigoplus_j (H^{i-j}(\mathbb P^n,\cal G(j))\otimes\Omega^{-j}(-j))$,
called {\em Beilinson's monad}, has cohomology $H^i(\cal K^*)$ equal to
$\cal G$ in degree $i=0$ and $0$ in all other degrees,
cf. \cite{Bei}.
\end{rem}

\bigskip
In the case of even sets of nodes on sextics, \cite{CaCa} classifies
the sets of cardinalities $t=24,32,40$.
These cases are given by the following symmetric bundle maps
(here and in the rest of the paper, we denote  the sheaf $\cal
O_{\PP^3}$ simply by $\cal O$):
$$
\begin{tabular}{ll}
$t=24$ & $\varphi: \cal O(-5)\oplus\cal O(-4) \rar \cal
O(-2)\oplus\cal O(-1)$\\
$t=32$ & $\varphi: 3 \cal O(-4) \rar 3 \cal O(-2)$\\
$t=40$ & $\varphi: \Omega^2(-1)\oplus \cal O(-4) \rar
\Omega^1(-1)\oplus \cal O(-2).$
\end{tabular}
$$

The cases of even sets of nodes with cardinality $t=48$ and $t=64$
are excluded by corollary (\ref{weights}).
Thus the only remaining case, in order to complete the classification
of even sets of nodes on sextic surfaces, is the case $t=56$ .

\bigskip
Let us then restrict ourselves now to the case
$d=6$ and $t= 56$, and let us consider the module
$M:= U \bigoplus\, (\oplus_{i>1} H^1(F, \cal F(i)))$.

We shall use  the analysis done in \cite[p. 254]{CaCa} for  
$\cal F\mid_H$,
where $H$ is a smooth plane section of $F$: it shows that
$h^0(F,\cal F(1))=h^2(F,\cal F(1))=0$
(in loc. cit. it is shown that this holds unless  $\cal F\mid_H$ is
of type (2,4),
but if $\cal F\mid_H$ of type (2,4) then $F$ is  of type (2,4) too,
in the sense of \cite[Thm. 2.2 and Thm. 2.16]{Ca},
and $t=24$).

Hence we can assume that $h^0(F,\cal F(1))=h^2(F,\cal F(1))=0$,
so that by Riemann-Roch applied to $\tilde F$
$$h^1(F,\cal F(1))=-\chi(\cal F(1)=-8+\frac{t}4=6.$$

According to the notation in \cite[p. 254]{CaCa},
set $2\tau:=h^1(F,\cal F(1))$, $a:=h^1(F,\cal F)=h^1(F,\cal F(2))$,
$b:=h^2(F,\cal F)=h^0(F,\cal F(2))$, 
the equalities following by Serre duality.
Our previous calculation yields $ \tau = 3$.

The exact sequence
\begin{multline}
         H^0(F,\cal F(1))=0 \rar H^0(H,F(1)\mid_H)
\rar H^1(F,\cal F) \rar H^1(F,\cal F(1)) \rar \\
         \rar H^1(H,F(1)\mid_H)\rar H^2(F,\cal F)
\rar H^2(F,\cal F(1))\cong H^0(F,\cal F(1))=0
\end{multline}
gives the relation $\chi(\cal F(1)\mid_H)-a+2\tau+b=0$.
An application of Riemann-Roch on $H$ yields
$\chi(\cal F(1)\mid_H)=-3$
and the above relation becomes $ 0 \leq b=a- 3$.

Finally notice that (see \cite[Formula 3.2, p. 248]{CaCa})
$$H^1(F,\cal F(-m))\cong H^1(F,\cal F(m+2))^\lor=0,\ m>0,$$
and trivially also $H^0(F,\cal F(-m))=0,\ m>0$.

Since the rank of $\cal F$ at the generic point of $\Pn 3$ is $0$,
a computation of the ranks of all the terms of the Beilinson's monad
of $\cal F(3)$ yields the relation
$4b+6\tau-4a-c=0$, i.e. $c=12-2\tau=6$.

Therefore, the Beilinson table of $\cal F(3)$ is:
$$
\vbox{\offinterlineskip
\halign{
\hbox to 10pt{\hfil}#
&\vrule height10.5pt depth 5.5pt#
&\hbox to 35pt{\hfil$#$\hfil}
&\vrule#
&\hbox to 35pt{\hfil$#$\hfil}
&\vrule#
&\hbox to 35pt{\hfil$#$\hfil}
&\vrule#
&\hbox to 40pt{\hfil$#$\hfil}
&\vrule#
&\hbox to 40pt{\hfil$#$\hfil}
&\vrule#
&#
\cr
\omit&\omit&&\omit
\hbox to 0pt{\hbox to 0pt{\hss$\big\uparrow$\kern
-0.4pt\hss}\raise5pt\hbox to 0pt{\ $i$\hss}}\cr
&\multispan{2}\hrulefill&\multispan{8}\hrulefill\cr
&& c=6 && b && 0 && 0 && 0 && \cr
&\multispan{2}\hrulefill&\multispan{8}\hrulefill\cr
&& 0 &&  b+3 && 2\tau = 6&& b+ 3 && 0 && \cr
&\multispan{2}\hrulefill&\multispan{8}\hrulefill\cr
&& 0 && 0 && 0 && b && c= 6 &&\cr
\multispan{12}\hrulefill&\kern -2pt
\hbox{\kern .4pt \vbox to 0pt{\kern 1.1pt\vss\hbox to
0pt{$\longrightarrow$ \hss}\vss}
             \vbox to 0pt{\kern 3pt \hbox to 0pt{\ $j$ \hss}\vss}}\cr
}\kern 6pt}
$$

\begin{prop}\label{h^0(2H-L)}
Let $F$ be a nodal surface of degree $6$ with $\overline t$ nodes and with
an even set $\Delta$ of
$t= 56$ nodes and
$L\in Pic(\tilde F)$ the corresponding divisor such that $\sum_{i\in
N_\Delta} A_i=2L$.
Then $b := h^0(2H-L)\leq 1$.
\end{prop}

\begin{proof}
Assume $h^0(2H-L)\geq 2$ and write $\mid 2H-L\mid=\mid M\mid + \Psi$, where
$\Psi$ is the fixed part of the linear system $\mid 2H-L\mid$.
Let $\mathcal L$ be the Gaussian linear subsystem of
$5H-\sum_{i=1}^{ \overline t} A_i$
and $C$ any effective divisor in $|2H-L|$.

Since $F$ is nodal it follows by our previous argument that $\L$
is free from base points, hence for any effective divisor $C'\leq C$ we have
$C' \L\leq C \L =  60-t=4$.

Observe that by \cite{JR} the number  $\overline t$ of nodes of $F$ satisfies
$\overline t\leq 65$ :
since $\mathcal L^2\geq 150-2\overline t\geq 20$,
the  Index Theorem ensures that
$(C')^2\leq 0$. In particular, it follows that $M^2=0$.

Since the dual surface $F^\lor$ is birational to $F$ and therefore it is still
of general type, $\gamma(M)$ has degree at least 4.

Hence, by the previous calculation, $M \L=4$ and $\Psi \L=0$,
the general curve in $|M|$ is irreducible (hence,  smooth) and maps
$1:1$ to a quartic.

Since the arithmetic genus of $M$ is at least 2, it follows that $M$
maps $1:1$ to a
plane quartic, and therefore its arithmetic genus is at most 3.

But then $4\geq 2p_a(M)-2=M \cdot (M+K_F)=M \cdot K_F=2M \cdot H$,
i.e. $M \cdot H\leq 2$, a
contradiction since $F$ is of general type.

\end{proof}

We can summarize the above results in the following statement.

\begin{thm}\label{HilbertFunctions}
Let $F$ be a nodal surface of degree $6$ with an even set of $56$ nodes.
Then $M$ is an Artinian module of length $2$ with Hilbert function
$( \tau, a)$ =
$(3,3)$ or $(3,4)$.
\end{thm}

\begin{proof}
By hypothesis $\tau=3$. Prop. (\ref{h^0(2H-L)}) yields $b=h^1(F,\cal
F(2))\leq1$.
Therefore only the following two cases are possible:
$b=0, a=3$ or $b=1, a=4$.
\end{proof}

We shall treat in the sequel only the first case, for the second we limit
ourselves to  parenthetically posing the following

\begin{prob}\label{problem}
Let $F$ be a normal sextic surface in $\Pn 3$ with an even set of $56$ nodes.
Can the  module $M:=U\bigoplus\,  (\oplus_{i> 1} H^1(F, \cal F(i)))$
have Hilbert function $(3,4)$ ?
\end{prob}

\begin{rem}
This case cannot be excluded by coding theory since there exists
a $9$-dimensional code
$K \subset ( \ZZ / 2 \ZZ)^{56}$ with weights $(24,32,56)$.
\end{rem}

\proof
This code is constructed as follows: consider a code $ U \subset (
\ZZ / 2 \ZZ)^{51}$
of dimension $8$ and weights $(24,32)$,
and let $e \in ( \ZZ / 2 \ZZ)^{56}$ be the vector with all
coordinates equal to $1$.
It suffices to define $ K$ as the span of $ U$ and $e$.

The existence of $U$ (cf.  \cite{mcw}, page  229) is easily
established if we let
$\FZ$ be the finite field with $ 2^8$ elements, $\xi$ a generator of $\FZ ^*$.

Then $\xi^5$ is a primitive $51$-st root of unity, it generates $\FZ$
as a field,
thus $\FZ\cong  (\ZZ / 2 \ZZ ) [\xi^5] / (P)$, where $P$ is an
irreducible polynomial
of degree $8$ dividing $ x^{51} - 1$. By the Chinese remainder theorem
       $\FZ\cong ( \ZZ / 2 \ZZ ) [\xi^5] / (P)$ is a direct summand of
$ ( \ZZ / 2 \ZZ) [x] / (x^{51} - 1) \cong  ( \ZZ / 2 \ZZ)^{51}$, and
it suffices to let
$U$ be the subspace of $ ( \ZZ / 2 \ZZ)^{51}$ which  corresponds to $\FZ$.

\qed

\section{ Hilbert function $(3,3)$: general features.}\

We shall assume, throughout the rest of the paper, that we have an even set
$\Delta$ of $56$ nodes, and that $b=0$, i.e.,
$a=3$  (cf. \ref{HilbertFunctions}) .
In other terms, the Artinian module $M$ has dimension $ \tau = 3$ in
degree $1$,
dimension $ a = 3$ in degree $2$, and $0$ in  degree $ \neq 1,2$.

\bigskip
By Theorem (\ref{CaCaThm}) (applied to $\cal F(3)$ instead of $\cal F$,
i.e., we replace the previous $\mathcal E$  by its twist $\mathcal E (3)$)
we have a resolution of $\cal F(3)$  of the form:
\begin{equation}
         0 \rar \mathcal E^{\vee} \rxar{\ \varphi} \mathcal E \rar
\mathcal F(3) \rar 0.
\end{equation}

In this setting, the symmetric map $\varphi$ appearing in the resolution of
$\cal F(3)$
belongs to $H^0(\mathbb P^3, S^2(\cal E))\subset Hom(\cal E^\lor,\cal E)$.

\begin{defn}
Throughout the rest of the paper we denote by $U$
a given Lagrangian 3-dimensional subspace  
of $H^1(\cal F(1)) $, and we denote by $W$ the    3-dimensional space
$W:=H^1(\cal F(2)).$

Moreover, we shall  
denote by  $V$ the four dimensional vector space
 $ V  : = H^0(\cal O(1))$.
Later on, more generally,  $V$ shall denote a four dimensional vector space
and we shall often continue to denote 
$Proj (V)$ by $\PP^3$.

Thus we shall have
\end{defn}

\begin{rem}
Beilinson's theorem and the cohomology table for $\cal F(2)$, 
which we described above, implies that
$\cal E(-1)$ is obtained by adding a direct sum of line bundles to
$$
\cal E'(-1):=\ker(U\otimes\Omega^1(1)\cong 3\Omega^1(1)\rar
W\otimes\cal O\cong 3\cal O),
$$
and that (since the Beilinson's complex has
 no cohomology in degree $\neq
0$) the above map
is surjective: hence 
$\cal E'
$ is a vector bundle with $rk(\cal E')=6$.

Consider now the Euler sequence
\begin{equation}
0\rar \Omega^1(1)\rar V\otimes\cal O\cong 4\cal O\rar\cal O(1)\rar 0.
\end{equation}
It implies that $h^0(\Omega^1(1))=0$ and $h^0(\Omega^1(2))=6$, thus
$h^0(\cal E'(-1))=0$ and, since
      Beilinson's Table for $\cal F(3)$ implies that $h^1(\cal E')=0$,
we infer that $h^0(\cal E')=3\times 6-4\times 3=6$.

On the other hand,  Beilinson's complex for $\cal F(3)$
yields an exact sequence:

$$0 \rar 3\cal O(-4) \rar 6\Omega^2(2) \rar 3\Omega^1(1) \oplus 6\cal O
      \ra \cal F(3)\rar 0,$$
and we make the following simplifying

\medskip\noindent
{\bf  FIRST ASSUMPTION:  $\cal F$ is generated in degree $3$}
and the linear map $H^0(\cal E') \ra H^0(\cal F(3)) $ 
is an isomorphism.

\end{rem}

\begin{prop}\label{1as}
According to the previous notation, 
the above first assumption implies that $\cal E=\cal E'$, equivalently, that
{\bf  rank   $(\cal E)$  = 6}. More precisely, it means that
there exists a homomorphism $\beta :  U\otimes\Omega^1(2) \cong
3\Omega^1(2) \ra  W\otimes\cal O(1) \cong 3\cal O(1) $ such that $\cal E =
ker
\beta$ and that we have an exact sequence
\begin{equation}\label{E=ker}
0\rar \cal E \rar U\otimes\Omega^1(2) \cong
3\Omega^1(2)  \rxar\beta W\otimes\cal O(1) \cong 3\cal O(1) \rar 0.
\end{equation}

Conversely, if $\cal E$ is obtained in this way, it is a rank 6 bundle
with an intermediate cohomology module $M$ with the required
Hilbert function
     of type $(3,3)$. Moreover  $H^0(\E^\vee) = 0$.
\end{prop}

\proof
If $\cal F$ is generated in degree $3$
there is an exact sequence
$$0\rar \tilde G \rar 6\cal O \rar \cal F(3) \rar 0.$$
where $h^0(\tilde G)= h^1(\tilde G)=0$ (cf. Beilinson's table).

Dualizing the sequence
$0\rar \cal E' \rar 3\Omega^1(2) \rxar\beta 3\cal O(1) \rar 0$
yields
$$  0\rar 3\cal O(-1) \rar 3 T(-2) \rar {\cal E'}^\lor \rar 0.$$
Thus $h^0( {\cal E'}^\lor)=0$. 

Assume now that $H^0_* (\cal E') \ra H^0_* (\cal F (3)) $ is not surjective.
Then, since by our assumption  
$H^0 (\cal E') \ra H^0 (\cal F (3)) $ is surjective,
$\cal E$  will be obtained from $\cal
E'$ by adding a direct sum of line bundles $\hol (-m)$ where $ m >0$.
This leads however to a contradiction,
 since then $\hol (m)$ is a direct summand
of $\cal E^\lor$ but  it cannot embed in 
$\E$ since $ H^0(\cal
E'(-1)) = 0$ implies that $ Hom (\hol (m), \cal E) = H^0(\cal
E(-m))  = H^0(\cal E'(-m)) = 0$.

\medskip
For the converse, we simply observe that if  (\ref{E=ker}) holds,
then $ H^1 ( \cal E (-2)) \cong 3 H^1 (\Omega^1),
H^1 ( \cal E (-1)) \cong 3 H^0 (\hol). $
Since $\E' \cong \E$ follows rightaway that $H^0(\E^\vee) = 0$.

\QED

Therefore we get the following exact commutative diagram:
$$\xymatrix{
& 0 & 0\ar[d]\\
0 \ar[r] &\cal E^\lor \ar[r]^\varphi \ar[u]&\cal E \ar[r] \ar[d]
        &\cal F(3) \rar 0\\
&U^\vee\otimes T(-2) \ar[u] \ar[r]^\Phi & U\otimes\Omega^1(2)
\ar[d]^\beta\\ &W^\vee\otimes\cal O(-1) \ar[u]_{^t\!\beta}  
&W\otimes\cal O(1)\ar[d]\\ & 0\ar[u] & 0\\}$$
and the map $\varphi$ yields, by composition, a homomorphism
$$\Phi\in Hom(U^\vee\otimes  T(-2), U\otimes \Omega^1(2))$$
which is symmetric since $\varphi$ is symmetric.
Conversely,
such a homomorphism $\Phi$ determines $\varphi$ if and only if
$\beta \Phi=\Phi\
^t\!\beta=0$: since however we  choose $\Phi$  symmetric  the two
conditions are
equivalent to each other.

A more concrete way to setup the parameter space for such vector bundles
      is to replace $Hom( T(-2),\Omega^1(2))$
via matrices of polynomials, as
follows.

Recall that   $V:=H^0(\Pn 3,\cal O(1))$ is the space of linear forms on $\Pn
3$. Applying $Hom(-,\cal O)$ to the Euler sequence and tensoring by
$\hol(1)$ yields, since $Hom(\hol(2),\hol(1))=0$,
$Ext^1(\hol(2),\hol(1))=0$,
\begin{equation}\label{B=beta}
Hom(\Omega^1(2),\cal
O(1))\cong Hom( V \otimes \cal O(1) , \cal O(1)).
\end{equation}

Thus the map $\beta$ factors through a map
$B: U\otimes(V \otimes\hol(1)) \rar W\otimes\hol(1)$ and
the sheaf map $B$ is surjective.
This surjectivity is obviously equivalent to 
$H^0(B(-1)):\ U\otimes V \ra W$ being surjective.
In the sequel we shall often identify the sheaf map $B$ with the
corresponding tensor $H^0(B(-1)) \in U ^{\vee}\otimes V^{\vee} \otimes
W$.

Let  $\epsilon$ be the tensor product of the identity map of the
isotropic subspace $U$ with
the evaluation map
$V\otimes \cal O\rar \cal O(1)$.
Then one sees easily that
$\cal E =\ker \beta=\ker B\cap\ker \epsilon$,
the short exact sequence (\ref{E=ker}) becomes
\begin{equation}\label{E=ker2}
      0 \rar \cal E \rar U\otimes V \otimes\cal O(1)) \rxar{(B,\epsilon)}
(W\otimes\cal O(1))\oplus (U\otimes\cal O(2)) \rar 0,
\end{equation}
     and the previous
diagram is replaced by:
$$\xymatrix{
& 0 & 0\ar[d]\\
0 \ar[r] &\cal E^\lor \ar[r]^\varphi \ar[u]&\cal E \ar[r] \ar[d]
&\cal F(3) \rar 0\\
& U^\vee\otimes V^\vee \otimes\cal O(-1)) \cong 12 \cal
O(-1) 
\ar[u] \ar[r]^A &U\otimes V \otimes\cal O(1) \cong 12 \cal
O(1)\ar[d]^{(B,\epsilon)}\\ 
&(W^\vee\otimes\cal O(-1))\oplus (U^\vee\otimes\cal O(-2))
\ar[u]_{^t\!(B,\epsilon)}
&(W\otimes\cal O(1))\oplus (U\otimes\cal O(2))\ar[d]\\
& 0\ar[u] & 0\\}$$
where by a similar token to the one before the map $\varphi$ yields
a symmetric matrix
$A\in Mat(12\times12,Hom(\cal O(-1),\cal O(1)))$ and conversely such
a matrix determines
$\varphi$ if and only if $(B,\epsilon)\cdot A=0$.

Thus we obtain, as a parameter space for the symmetric resolutions of
$\cal F(3)$
satisfying the open condition  given by the main assumption,
the variety of pairs
\begin{equation}\label{AB}
\begin{split}
\frak M_{AB} : = &\{(B,A)\mid
B\in Mat(3\times 12,\CZ),\\&
A\in Mat(12\times 12,H^0(\cal O_\Pn 3(2))),
A=^t\!A, (B,\epsilon) \cdot A=0\}.
\end{split}
\end{equation}

As a matter of fact, the equation of the surface $F$ will then be given 
as the G.C.D. of the determinants of the $ 6 \times 6$ minors of the matrix $A$, 
whereas the even set of nodes $\Delta$ will be found to be given by 
the ideal of the determinants of the $ 5 \times 5$ minors of the matrix $A$.

A direct but complicated calculation
shows that for a general choice of the parameter $B$ (determining the
bundle $\cal E$)
the solution space  for the $A$'s (yielding the symmetric map
$\varphi$) has positive
dimension.

However, by computer algebra, one checks  that  a
random choice of $B$
and a random choice of $A$ do not give a sextic surface with 56 nodes,
but the square of a cubic surface.
Notice that the condition that a given pair ($B,A$) yields a sextic with 56
nodes is an open condition, and therefore
there can exist sextics with an even
set of 56 nodes and satisfying the main assumption only if
the above parameter space
is reducible: we shall later show that this is indeed the case.

We finish this section by remarking that $B$ is the
multiplication matrix of the module $H^1_*(\cal E)$
     (i.e., the matrix of the only  part of the multiplication map which is not
a priori trivial).

\begin{rem}
The cohomology  exact sequence associated to the following 
twist of (\ref{E=ker2}), namely:
$$ 0 \rar \cal E(-2) \rar U\otimes V\otimes \hol(-1) \rxar{(B,\epsilon)}
      (W\otimes \hol(-1)) \oplus (U\otimes \hol) \rar 0$$
yields a canonical isomorphism $ U\cong H^1 (\cal E(-2))$.

Since there is a canonical isomorphism
$$H^0(\epsilon(-1)): U\otimes V \ra U\otimes H^0(\hol(1)),$$
the projection of
$W\oplus(U\otimes V) \ra W$ induces an isomorphism of
$H^1(\cal E(-1))= Coker H^0((B \oplus \epsilon)(-1))$ with 
$W$, such that the map 
$B:\ U\otimes V \rar W$
corresponds to the multiplication map of the
cohomology module $ H^1_{*} (\E)$. 
\end{rem}

\begin{rem}\label{b}
The condition that the linear map $B$ has maximal rank = 3 (which, as we
observed, follows from the first assumption) is obviously equivalent to
the condition that the module $M$ is generated in degree $-2$.
On the other hand, it also implies that there is an exact sequence

$$ 0 \ra \cal E \ra 9 \hol (1) \rxar b 3 \hol (2) \ra 0.$$
\end{rem}

\begin{rem}\label{conj1}
Let $F$ be a nodal sextic surface in $\Pn 3$ with an even set of $56$ nodes,
and assume that  $h^1(F, \cal F (1)) = 6$.
One may ask whether there is a
Lagrangian subspace $U$ such that the module
$ M =  U \oplus W :=  U \oplus H^1(F, \cal
F (2))$, which  has Hilbert function
$(3,3)$,  is generated by $U$.

A necessary condition is that, setting $ \tilde{W}: = H^1(F, \cal F (1)) $,
$ \tilde{W} \otimes V \ra W$ be surjective. In turn this is equivalent
to the pairing $ H^1(\hol_S ) \times H^0(f^* \hol_F (1)) \ra
H^1(f^* \hol_F (1))$ being non degenerate in the first factor.
\end{rem}

We proceed in the next section with the analysis of the vector bundles
corresponding to a general choice of $B$, giving a
geometrical explanation of the phenomenon
of which the computer made us aware.

\section{General Bundles and cubic surfaces}

Main purpose of this section is  to describe the beautiful geometry
which relates the main component of the moduli space of our vector
bundles with given intermediate cohomology module $M$ and the space
of cubic surfaces viewed as blow ups of the projective plane in six points.

Let us preliminarily observe that, if the main assumption is satisfied,
the vector bundle $\E$ is determined by the matrix $B$, hence we have
an irreducible parameter space for our vector bundles, and each open condition,
if verified at some point, is verified by the generic bundle $\E$.

Next, we have a surjection $H^0(\E) \ra H^0(\F(3))$ and we have seen
that both spaces are 6-dimensional, whence we get a homomorphism
$\i : 6\cal O \rar \cal E $. We make the

\medskip\noindent
{\bf SECOND ASSUMPTION:\\
1) $ \i : 6\cal O \rar \cal E $ is injective},
whence  an exact sequence:
\begin{equation}\label{seq1a}
           0\rar 6\cal O \rar \cal E \rar \tau \rar 0,
         \end{equation}
{\bf 2)
the torsion sheaf $\tau$  is $\cal O_G$-invertible,} where $G$ is the
divisor of $\Lambda^6(\i)$.

\medskip

\begin{lem}\label{ChernC}
Let a vector bundle  $\E$  be given as in \ref{E=ker}
or as in \ref{E=ker2}
Then its total Chern class is
\begin{equation}\label{c(t)}
     c(\E) (t) = 1 + 3 t + 6 t^2 + 4 t^3.
\end{equation}
In particular,  if the second assumption is
satisfied, {\bf  the divisor $G$ is a cubic
surface}.
\end{lem}

\proof
The sheaf $\tau$ has
      Chern polynomial 
$c(\tau)=c(\cal E)=c(\Omega^1(2))^3 c(\cal O(1))^{-3}=
(c(\cal O(1))^4 c(\cal O(2))^{-1})^3c(\cal O(1))^{-3}=
(1+t)^9(1+2t)^{-3}=(1+9t+36t^2+84t^3)(1-6t+24t^2-80t^3)=
1+3t+6t^2+4t^3$.

\qed

\begin{rem}\label{Rem:DoubleCubic}
Observe that the space $H^0 (S^2(\E))$ of symmetric morphisms from
$\cal \E^\lor$ to $\cal E$
contains $H^0(S^2(6\cal O))$ since to
$\tilde\alpha\in H^0 (S^2(6\cal O))$ corresponds $\alpha : =\iota^\lor \tilde
\alpha \iota$.
For these morphisms one has $\det\alpha=\det\tilde\alpha\, (\det \iota)^2$,
whence in this case $ div (det (\alpha)) = 2 G$,
and not  a  sextic surface.
\end{rem}

The next lemmas are meant to investigate the question: when does one have
equality  $h^0 (S^2(\E)) = 21$, i.e., when is $H^0 (S^2(\E)) =
H^0(S^2(6\cal O))$?
\smallskip

In order to answer this question, it is convenient first to analyse
the geometry and
the cohomology of the invertible sheaf $\tau$ on $G$.

\begin{rem}\label{tau}
Even without assuming $\tau$ to be $\hol_G$-invertible,
set $\tau'=\cal Ext^1(\tau,\cal O)$:
then the dual of the previous exact sequence
\ref{seq1a} gives
         \begin{equation}\label{seq1b}
0\rar \cal E^\lor \rar 6\cal O \rar \tau' \rar 0.
         \end{equation}
and we have:
\begin{enumerate}
\item By (\ref{seq1a}) clearly holds $H^0(\tau)=H^1(\tau)=H^2(\tau)=0$.

\item From (\ref{seq1b}) and $h^i(\cal E^\lor)\cong h^{3-i}(\cal E(-4))$ we get
        $h^0(\tau')=6$, $H^1(\tau')=H^2(\tau')=0$.
\item Since by definition $\tau'=\cal Ext^1(\tau, \cal O)$, applying
the functor
$\cal Hom(\tau,-)$ to the exact sequence
$0\rar \cal O\rar \cal O(3) \rar \cal O_G(3)\rar 0$ we get
$\tau'=\cal Hom (\tau,\cal O_G(3))$.
Therefore, if $\tau=\cal O_G(D)$, then $\tau'=\cal O_G(3H-D)$.

\end{enumerate}
\end{rem}
Since  $h^i(D)=0 \forall i$,
$h^0(3H-D)=6$, $h^i(3H-D)=0$
for $i=1,2$, follows by Riemann Roch that
$D^2+DH=-2$ and $10=36-7DH+D^2$. Therefore $HD=3$, $D^2=-5$.

By setting $\Delta:=D+H$, it turns out that $\Delta
H=6,\Delta^2+\Delta K_G=-2$.

\begin{lem}\label{D}
Assume that $G$ is a smooth cubic surface: then
there exists a realization of  $G$ as the image of the plane 
under the system
$| 3L - \sum_1^6 E_i|$ of plane cubics through six points, such that either
$\Delta \equiv 2L$, i.e., $\Delta$
corresponds to the conics in the plane, or (up to permutations of the six
points) 
$\Delta \equiv 3L - 2 E_1 - E_2.$
\end{lem}

\Proof
Observe preliminarly that if $ |H| = | 3L - \sum_1^6 E_i|$
is such a planar realization of a cubic surface, then another one
is obtained via a standard Cremona transformations
centered at three of the points $P_i$ corresponding to the 
$-1$-curves $E_i$.

In fact, if $ L' : = 2L - E_1 - E_2 - E_3 $, then 
$$ H =  3L'  - (L - E_1 - E_2) -  (L - E_1 - E_3) - (L - E_2 - E_3)
- E_4 - E_5 - E_6. $$

We have $0= H^2(D)= H^0(-D-H)$ and
 a fortiori $H^2(\Delta)=H^0(-D-2H)=0$.
It follows that
$|\Delta|$ has $h^0(\Delta)\geq 6$, $\Delta^2=4$, and the arithmetic
genus $ p_a(\Delta)=0$.

Since $ 0 = H^1(D)= H^1( K_G-D) = H^1(-H -D) = H^1(- \Delta)$, 
 it follows that $\Delta$ is connected. 

Hence we have a representation
$\Delta \equiv n L - \sum_1^6  a_i E_i$, where the $a_i$'s are non negative
and we assume $a_1\geq a_2\geq \dots\geq a_6$.

We have: $ \Delta^2 = 4 = n^2 - \sum_1^6  a_i^2  $, $\Delta
\cdot H=6 = 3 n - \sum_1^6  a_i$, i.e.
      \begin{equation}\label{nn2}
n^2 = \sum_1^6  a_i^2 + 4 , \  3 n = \sum_1^6  a_i + 6 .        
      \end{equation}

We want to show that, after a suitable 
sequence of standard Cremona transformations,
$\Delta\equiv 2L$  or $\Delta\equiv 3L-2E_1-E_2$.
By \ref{nn2}, we have $n\geq 2$ 
and for $n=2,3$ $\Delta$ is as claimed.
Hence the claim is that there exists a 
sequence of standard Cremona transformations
which makes  $|\Delta|$
have degree $n\leq 3$.

By applying $|2L-E_1-E_2-E_3|$ we get a new system $\Delta'$
with degree $n'=2n-a_1-a_2-a_3$.

By our ordering choice for the $a_i$'s, we have 
$$a_1+a_2+a_3 \geq (\sum_{i=1}^6 a_i)/2=3n/2 -3,$$
with strict inequality unless all $a_i$'s are equal.
We study this latter case  first:

\bigskip
{\bf Sublemma.} In the previous setting, $a_1=a_2=\ldots=a_6$ if and only
if $n=2$ and $a_1=a_2=\ldots=a_6=0$ or $n=10$ and $a_1=a_2=\ldots=a_6=4$.

\medskip
\proof
The statement follows immediately by defining $a:=a_1=a_2=\ldots=a_6$ and
using both conditions of \ref{nn2}: $n=2a+2$, $n^2=6a^2+4$,
which imply $ 8 a =2 a^2$.
 
\qed

\bigskip

The previous inequality gives:
$$n'\leq\frac{n}2+3\leq n \text{\qquad for $n\geq 6$},$$ 
and $n'<n$ for $n\geq 6$ unless $n'=n=6$ and $a_1=a_2=\ldots=a_6$,
which has no solution by the above sublemma.
We conclude that after  suitable Cremona transformations $n\leq 5$.

If $n=5$, then $n'\leq 5/2+3$, i.e., $n'\leq 5$.
Moreover, if also $n'=5$, then $a_1+a_2+a_3=5$ and using again
\ref{nn2} we obtain $a_4+a_5+a_6=4$. But then $a_1=a_6+1$ 
and we easily get
a contradiction since then $a_2= a_3=a_4=a_5 $
and they either equal $a_1$ or $a_6$.
Hence, after a suitable Cremona transformation,
we can always reduce to the case $n\leq 4$.

Let now $n=4$. Using \ref{nn2} we get
$\sum_{i=1}^6 a_i=6$ and $\sum_{i=1}^6 a_i(a_i-1)=6$.
We have the following two possibilities:
$a_1=3, a_2=a_3=a_4=1, a_5=a_6 =0$ or $a_1=a_2=a_3=2,
a_4=a_5= a_6 =0$.
In both cases we have that $a_1+a_2+a_3\geq 5$, and therefore $n'\leq 3$.

%

\QED

\begin{rem}
The complete linear system $\Delta$ has as image in $\PP^5$ either
the Veronese embedding of $\PP^2$, or the embedding
of $\PP^1 \times \PP^1$ through  
$ H^0 ( \hol_{ \PP^1 \times \PP^1 } (1,2))$.
In both case we have a surface of minimal degree (=4).
\end{rem}

\bigskip
Thus we have concluded that either $D=2L-H=-L+\sum E_i$, or
$D=3L- 2 E_1 - E_2 - H=- E_1 + \sum_3^6 E_i$.

\begin{cor}
$H^i(\hol_G( 2D) =0$ for $ i=0,2$, $h^1(\hol_G( 2D) =6$.
\end{cor}

\begin{proof}
The second part follows from the first by Riemann Roch, for
the first it suffices to intersect with $L$ (using Serre duality
in the case of $H^2(\hol_G( 2D) $).
\end{proof}

We are now ready to show that the smoothness assumption
for the cubic surface $G$ implies that
all symmetric morphisms from $\cal E^\lor$ to $\cal E$ factor
through and are
induced by symmetric morphisms form $6\cal O$ to $6\cal O$, whence  their
determinant is a double cubic, instead of a nodal sextic
(cf. Rem. \ref{Rem:DoubleCubic}).

\begin{lem}\label{sym}
Let
$$0\rar \cal F \rxar{i} \cal E \rar \tau \rar 0$$
be a locally free resolution of a coherent torsion sheaf $\tau$, which is
$\cal O_G$-invertible on a divisor $G$.

Then we have an exact sequence
\begin{equation}\label{eq:S2E}
         0 \rar \Lambda^2\cal F\rar \cal F\otimes \E \rar S^2 \cal E \rar
\tau\otimes\tau \rar 0
\end{equation}
and a monad
\begin{equation}
         0 \rar S^2 \cal F \rar \cal F\otimes \cal E \rar \Lambda^2 
\cal E \rar 0
\end{equation}
whose cohomology in the middle is exactly $Tor^1(\tau,\tau)$.
\end{lem}

\begin{proof}

Recall that locally, by our assumption, we can write:
$$\cal E=\cal Oe_1\oplus \cal Oe_2\oplus \ldots \oplus \cal Oe_r,
\ \cal F=\cal Oxe_1\oplus \cal Oe_2\oplus \ldots \oplus \cal Oe_r,$$
where $x$ is a local equation for $G$.

\medskip
Since $Tor^1(\cal B,\cal B')=0$ if $\cal B$ is locally free,
we obtain the following  commutative diagram with exact rows and columns:

$$\xymatrix{
&&0&0&0&\\
0 \ar[r] &Tor^1(\tau,\tau) \ar[r] &\cal F\otimes\tau \ar[r] \ar[u]
&\cal E\otimes\tau \ar[r] \ar[u] &\tau\otimes\tau \ar[r] \ar[u]&0\\
& 0 \ar[r] &\cal F\otimes\cal E \ar[r] \ar[u] &\cal E\otimes\cal E
\ar[r] \ar[u] &\tau\otimes\cal E \ar[r] \ar[u] &0\\
& 0 \ar[r] &\cal F\otimes\cal F \ar[r] \ar[u] &\cal E\otimes\cal F
\ar[r] \ar[u] &\tau\otimes\cal F \ar[r] \ar[u] &0\\
& & 0 \ar[u]&0\ar[u] &Tor^1(\tau,\tau)\ar[u] &\\
&&&&0\ar[u] }$$

Hence the composite map
$\cal E\otimes\cal E \rar  \tau\otimes \cal E \ra \tau\otimes\tau$
is surjective and has kernel
generated by $(\cal F\otimes \cal E ) \oplus ( \cal E\otimes\cal F)$, as
the diagram shows.
Let $K$ denote the kernel of the  map
$(\cal F\otimes \cal E ) \oplus ( \cal E\otimes\cal F) \rar \cal
E\otimes\cal E$.
$K$ contains the image of $\cal F\otimes \cal F$ through the inective map
$(id\otimes i,-id\otimes i)$, where $i:\ \cal F \rar \cal E$ is the
inclusion.

Therefore we get the complex
$$
0\rar \cal F\otimes\cal F \rar (\cal F\otimes \cal E) \oplus (\cal E\otimes\cal
F)
\rar \cal E\otimes\cal E \rar \tau\otimes\tau\rar 0,
$$
exact except possibly at $(\cal F\otimes \cal E) \oplus (\cal E\otimes\cal F)$,
where the cohomology is equal to $K/(\cal F\otimes \cal F)$.

Now let $K_1$ be the inverse image  of $Tor^1(\tau,\tau)$ in $\cal
F\otimes \cal E$
via the short exact sequence
$0 \rar \cal F\otimes\cal F \rar \cal F\otimes \cal E \rar
\cal F \otimes \tau  \rar 0$.

We claim that $K\cong K_1$.
In fact, if $h_1\oplus(-h_2)\in K$, the diagram shows that $h_1\in K_1$ 
and moreover that $h_2$ is uniquely determined.
Conversely, if $h_1\in K_1$, then there is an (unique) 
element $h_2\in\E\otimes \F$ with $h_1=h_2\in \E\otimes \E$.

This implies that $K/{\F\otimes \F}\cong Tor^1(\tau,\tau)$.

Locally, $K_1$ is generated by $x e_1\otimes e_1  
(\mod{\F\otimes \F})$.
Thus $K/{\F\otimes \F}$ is generated by 
$(x e_1\otimes e_1)\oplus(-e_1\otimes xe_1)\in (\F\otimes \E)\oplus (\E\otimes \F)$,
which is an antisymmetric tensor, and we are done.

\end{proof}

\begin{cor}\label{21}
According to the previous notation,
assume that $G$ is a smooth cubic surface.
Then $h^0(S^2 \cal E)=21$, $h^1(S^2 \cal E)=6$, $h^2(S^2 \cal
E)=h^3(S^2 \cal E)=0$.

The same conclusion $h^0(S^2 \cal E)=21$ holds if  more generally
$\cal E$ verifies the first and second assumption
and $ H^0 ( \tau^{\otimes 2}) = 0$.
\end{cor}

\begin{proof}
Split the long exact sequence (\ref{eq:S2E}) into
$$0 \rar 15\cal O \rar 6\cal E \rar \cal H \rar 0,$$
$$0 \rar \cal H \rar S^2\cal E \rar \tau^{\otimes 2} \rar 0.$$
Recall from the construction of $\cal E$
that $h^0(\cal E)=6$ and $h^i(\cal E)=0$ for $i=1,2,3$.
       The first corresponding long exact   cohomology sequence yields
$h^0(\cal H)=36-15=21$
and $h^i(\cal H)=0$ for $i=1,2,3$.

Therefore it suffices to observe that if $G$ is smooth, by the
previous corollary,
one has $h^0(\tau^{\otimes 2})=h^2(\tau^{\otimes 2})=0$,
      $h^1(\tau^{\otimes 2})=6$. The rest is straightforward.

\end{proof}

Recall now that the vector bundle $\E$, provided that the main assumption
and the second assumption are satisfied, produces an invertible sheaf $\tau$ on
a cubic surface $G$; conversely, given such a sheaf $\tau$, one can
construct
     $\cal E$ as an extension of $6 \cal O$ and
$\tau$ as in \ref{seq1a}.

     Setting as before $\tau' : = \E xt^1(\tau,\cal O)$, we see that
such an extension is parametrized by
$Ext^1(\tau,6 \cal O)=H^0(6\cal Ext^1(\tau,\cal O))\cong \CZ^{36}$,
if, as in remark \ref{tau}, we have $h^0(\tau')=6$.

\begin{lem}\label{unicity}
Assume that $h^0 (\E^\vee ) = 0$ (cf. the proof of proposition 3.2),
and that $\E$ is an
extension  as in  \ref{seq1a}: then the extension class in
$Ext^1(\tau,6 \cal O)=H^0(6\cal Ext^1(\tau,\cal O))\cong \C^{6}
\otimes \C^6$ is a rank 6 tensor
( we shall refer to this statement by saying that {\bf the extension
does not partially
split}).

In particular,   {\bf $\E$ is then
     uniquely determined up to  isomorphism}.

\end{lem}
\proof

The extensions which yield vector bundles
form an open set.

     We canonically view these extension classes
as $ Hom(H^0(\tau'),H^0( 6 \cal O)) = Hom(H^0(\tau'),\C^ 6 ) $
through the coboundary map of the corresponding exact sequence.
We have then an action of $GL(6, \C)$ as a group of automorphisms of
$ 6 \cal O$ , which induces an action on $ Hom(H^0(\tau'),H^0( 6 \cal O)) =
Hom(H^0(\tau'),\C^ 6 ) $ which is immediately identified to
the composition of the corresponding linear maps.

The extensions which yield vector bundles
form an open set, which contains an open dense orbit,
on which this action is free, namely, the tensors of rank = 6.

If the rank of the tensor corresponding to an extension is $ = r < 6$,
it follows that the extension is obtained from an extension
$ 0 \ra r \hol \ra \E'' \ra \tau \ra 0 $ taking then a direct sum
with $ ( 6 - r) \hol$: but then $ ( 6 - r) \hol$  is a direct summand
of $\E^\vee$, a contradiction.

\qed

\begin{cor}\label{vb}
$\E$ as above (\ref{unicity}) is a vector bundle if
$H^0(\tau')$ has no base points.
\end{cor}

\proof
Our hypothesis shows that in each point of $G$ the local extension class is
non zero, hence it yields a locally free sheaf.
\qed

Let us now show that the second case contemplated in Lemma \ref{D}
does not occur, since it produces a vector bundle $\E$ with a different
intermediate cohomology as the one we require.

\begin{lem}
The second case in \ref{D} does not occur, else the associated
vector bundle
$\E$ would then have $h^2(\E (-3)) = 1 \neq b=0$.

\end{lem}
\Proof
Assume that $ D = \sum_3^6 E_i - E_1$: then the linear system
$ | 2 H - D | = | 6 L  - E_1 - 2 E_2 - 3 (\sum_3^6 E_i ) |$ has bigger
dimension than the expected dimension $ 27- 28 = -1$, since it contains
an effective divisor, $| 2 L  - E_1 -  (\sum_3^6 E_i ) | + 2 | 2 L  - E_2 -
(\sum_3^6 E_i ) | $. This amounts to the nonvanishing of the
cohomology group $H^2 ( - 3 H + D) = H^2 (\tau (-3))$.

    From the exact sequence \ref{seq1a} we infer that $h^2(\E (-3)) =1$,
whereas we assumed throughout that $h^2(\E (-3)): =b = 0$,
a contradiction.

\qed

We assume now  that $G$ is a smooth cubic surface,
and that $\tau$ is an invertible sheaf on $G$, corresponding to
the divisor class $ -L + \sum_1^6 E_i$.

Consider now the associated vector bundle $\E$ : we want to verify that
$\E$ has the required cohomology table, i.e., we want to calculate
the dimensions  $ h^i (\E (-n))$  for $ n=0,1,2,3$.
This will   allow us  to verify that there are bundles
$\E$ which satisfy the main and the second assumption.

\begin{lem}\label{expected}
Let $G$ be a smooth cubic surface,
and let $\tau$ be the invertible sheaf on $G$ corresponding to
the divisor class $ -L + \sum_1^6 E_i$: then
the associated vector bundle $\E$
has the required cohomology table.
\end{lem}

\proof
Observe that

\begin{itemize}

\item
clearly $h^0(\tau (-n)) = h^0( D - n H) = h^0( -(1 + 3n) L + (n-1)
\sum_1^6 E_i) = 0$,
     for $ n=0,1,2,3$
\item
clearly $h^2(\tau(-n))= h^0( (n-1) H - D ) =h^0( (3n -2) L - n
\sum_1^6 E_i)= 0$,
     for $ n=0,1,2,$
since a quartic with $6$ double points is a union of $4$ lines
\item
$h^1(\tau (-3))= h^1( - 3 H +  D ) =h^1( - (10 L - 4 \sum_1^6 E_i)= 0$
by Ramanujam's vanishing trick for regular surfaces,
since the linear system $ | 10 L - 4
\sum_1^6 E_i |$ contains a reduced and connected divisor,
namely $Q_1+\ldots+Q_5+E_6$,
where $Q_i\in|2L-\sum_1^6 E_j+E_i|$.
\item
     since $\chi(\tau(-n))=1+\frac12
(D-nH)(D-(n-1)H)= \frac32 n (n-3)$, we have also
$h^1(\tau)=0$,  $h^1(\tau(-1))=h^1(\tau(-2))=3$, $ h^2(\tau(-3))= 0$.
\end{itemize}

\QED

We have seen how  to a linear map $ B :  V \otimes U \ra W $,
where $V=H^0(\Pn 3, \cal O_\Pn 3(1))$ denotes the space of linear
forms, $U $ is a fixed isotropic subspace in $H^1 (\F (1))$,
$ W : = H^1 (\F (2))$, corresponds a homomorphism of
vector bundles
$ B :    V \otimes U \otimes \hol \ra W \otimes \hol$,
inducing $\beta : U \otimes \Omega^1 (2) \ra  W \otimes \hol(1)$,
whence finally a vector bundle $\E : = ker (\beta)$ if $\beta$ is
surjective.

The second assumption yields a cubic surface $G \subset \PP roj(V)$
and an invertible sheaf $\tau$ on $G$. If $G$ is smooth, 
the invertible sheaf
$ \tau (1)$ yields then a birational morphism onto
a Veronese surface, whence represents $G$ as the blow up of a projective
plane $\PP^2$  in a subscheme $\zeta$ consisting of six points, and
as the image of
$\PP^2$ through the linear system of cubic curves  passing through
$\zeta$. The Hilbert-Burch theorem allows us to make an explicit construction
which goes in the opposite direction.

\begin{rem}
Let $U',W'$ be 3-dimensional vector spaces.

Consider a $3 \times 3 \times 4$ tensor $\hat{B}\in(U'\otimes V\otimes
{W'}^\vee)$
and assume that the corresponding  sheaf homomorphism
$\hat{\cal B}  : W'\otimes \cal O_{\PP^2}(-1) \ra V\otimes \cal O_{\PP^2} $
     yields an
exact sequence
\begin{equation}
         0 \rar W'\otimes \cal O_{\PP^2}(-1) \rxar{\hat{\cal B}} V\otimes \cal O_{\PP^2}
\rar \cal O_{\PP^2}(3) \rar \cal O_\zeta(3)
\rar 0
\end{equation}
which is the Hilbert Burch resolution of a codimension $2$ subscheme 
$\zeta$
of length $6$.

We obtain a canonical isomorphism $V\cong H^0( I_\zeta(3))$ and we 
let  $G\subset\PP roj(V)$
be the image of $\PP^2$ via the rational map $\psi$ associated to $V$.
     Under the above
 assumption on  $\hat{B}$, and if moreover $\zeta$ is a local complete
intersection,
$G$ is a normal cubic surface, and if we define $\cal G:=(\psi_* (\cal O(1))$,
then we have an  exact sequence on
     $\PP roj(V)$:

\begin{equation}
         \label{Geq1}
         0 \rar W'\otimes \cal O (-1)  \rxar{\cal B} U'\otimes \cal O   \rar
\cal G \rar 0.
\end{equation}

\end{rem}
\proof
Let us continue to use our previous notation: $\hol_G(L) = \cal G$,
and $\hol_G(H) = \hol_G(1) $, so that our explicit choice of $\hat{B}$
provides a canonical basis $x_0,x_1,x_2$ of $H^0(\cal G)$,
and a canonical basis $y_0,\ldots,y_3$ of
$V = H^0(I_\zeta(3))$ provided by the cubic polynomials $\gamma_i (x)$
which are the coordinates of  $ \Lambda^3 (\hat{\cal B} )$;
in other words,
$(\gamma_0(x),\ldots,\gamma_3(x))$
     are the $3\times 3$ minors of $\hat{\cal B}$, taking with alternate signs,
and the Hilbert Burch sequence amounts to giving
     the relations
$$\sum_{k,i} (\hat{B})_{k,j}^i \ x_i \gamma_k (x)=0.$$
In turn, we view these as relations on $G$,
$\sum_{k,i} (\hat{B})_{k,j}^i \ x_i y_k =0.$

We define now the  $ 3 \times 3 $ matrix ${\cal B} $ as:
\begin{equation}\label{Btidle}
{\cal B}_{i,j} : =\sum_{k} (\hat{B})_{k,j}^i \  y_k.
\end{equation}
Let us  still denote by $\cal B$ the morphism in $Hom(W'\otimes
V^\vee,U')$ represented by the above matrix with respect to the chosen
bases:
$\cal B$ is another way of representing the $3 \times 3 \times 4$ tensor
$\hat{B}\in(U'\otimes V\otimes {W'}^\vee)$ and induces
then the exact sequence \ref{Geq1} on
     $\PP roj(V)$.

\qed

Under the more general assumption that $\cal B$
never drops rank by 2, $\cal G$ is
an invertible
sheaf on a cubic surface $G$, and $\cal G=\cal O_G(L)$, with $h^0(L)=3$.

Conversely,  given an exact sequence as \ref{Geq1},
the space $H^0(\cal G)$, since $\cal G$ is generated by global
sections, yields a
morphism
$\pi : G \ra \PP roj (H^0(\cal G))$.
We calculate the  Euler characteristic
of $\cal G$:
$$\chi(\cal G(n))=3\left[\binom{n+3}3-\binom{n+2}3\right]=\frac12
(3n^2+9n+6).$$
If the cubic surface $G$ is normal, and we let $G'$ be its minimal
resolution, the
Hilbert Polynomial of
$\cal O_{G'}(L)$ is equal to
$\chi(L+nH)=1 + \frac12 (L+nH)(L+(n+1)H)=
1+\frac12(L^2+LH) +\frac32 n^2 + n(\frac32+LH)$
and therefore $LH=3$, $L^2=1$.
     Since $L^2=1$ $\pi$ is a birational morphism, and since $L$ has genus $0$
and degree three, $G'$ is the  birational image of $\PP^2$ under a
linear system of
plane cubics $H^0 (\cal I _{\zeta} (3))$, where $\zeta$ is a length six
zero-dimensional subscheme.

\begin{prop}\label{gen}
More generally, assume that $\cal G$ is an invertible sheaf on a
(non necessarily irreducible) cubic surface $G$, given by \ref{Geq1}: then
the Cartier divisor $L$ has degree $3$, and  the morphism $\pi$
yields a birational
morphism of a component of the
surface $G$ onto the plane.

Assume that $G$ is irreducible and that $ \tau : = \cal G^{\otimes 2} (-1)$:
then $ H^0 (\tau ^{\otimes 2} )= 0.$
\end{prop}

\Proof
It suffices to consider two general divisors $L_1, L_2$ in the linear
system $|L|$,
and to consider the Hilbert polynomials of the sheaves appearing in
the two exact sequences
$$0 \ra \hol_G(nH) \ra \hol_G(nH+L_1) \ra \hol_{L_1}(nH + L_1) \ra 0$$
$$0 \ra \hol_{L_2}(nH) \ra \hol_{L_2}(nH+L_1) \ra \hol_{L_1\cap L_2}(nH + L_1)
\ra 0.$$
The conclusion is that  $H^0 (\hol_{L_1\cap L_2}(nH + L_1)) = 1
\forall n >> 0$,
thus $L_1 $ and $L_2$ meet transversally in a smooth point.

Let us now assume that $G$ is irreducible, and consider the inverse of
the birational morphism $\pi$. We can factor it as a sequence
of blow ups $\sigma : Y \ra \PP^2$ followed by a  projection
$ p : Y  \ra G$, so to it corresponds a sublinear system of
a complete linear system on $Y$, which reads out on the plane
as $|H| = | b L - \sum_i a_i E_i|$. Here, $  b = H \cdot L$, and if
$ b \geq 3$, clearly $ L \cdot ( 4 L - 2 H)  <  0$, hence
$ | 4 L - 2 H | = \emptyset$ and our desired vanishing is proven.

If instead $ b \leq 2$, since $ dim |H| \geq 3$, it follows
that $ H = 2 L - E$ and $G$ is a linear projection
of the cubic scroll $ Y \subset \PP^4$ with centre a point
in $ \PP^4 \setminus Y$.
We claim however that this case  does not occur.

Essentially, since otherwise $ \cal G = p_* ( \hol_Y (L))$
would not be invertible.

As an alternative argument, observe that the factorization
$\sigma = \pi \circ p$ is not possible,
since $\sigma$ is an isomorphism on the complement of the
line $ E \subset Y$, while the inverse image of the double curve
of $G$ is a conic (possibly reducible) contained in $Y$.

\QED

\begin{rem}\label{doubleline}
The case where $G$ is a a linear projection of the cubic scroll
yields two sheaves $ \cal G$ which are not invertible.

As it is well known, every point in $\PP^4$ lies in one of the planes
spanned by the conics of the system $L$. 
If we project from a point in $\PP^4\setminus Y$, this conic maps two to
one to the double line of the cubic $G$.

Such a plane is said to be {\bf special} if the conic splits
into two lines $ E + F, F \equiv L-E$.

In the non special case,  we may assume without loss of generality
that the conic corresponds to the line $ z = 0$ in the plane,
that the blown up point is the point $ x=y=0$, and that
the linear system mapping to $G$ is generated by
$ ( zx : =  x_0,   zy :  = x_1, x^2 : = x_2, y^2 : = x_3) $.
In this case one sees that the matrix $\cal B$ is
$$\cal B=\begin{pmatrix}
x_0 &0 &x_1\\ 0 &x_1 &-x_0\\ -x_2 &-x_3 &0
\end{pmatrix},$$
     and that $G$ is then the cubic of equation
$-x_0^2x_3+x_2x_1^2=0$.

In the special case, we may again assume that
the blown up point is the point $ x=y=0$, we assume that the line $F$ is
the proper transform of $x=0$,
and that the linear system mapping to $G$ is generated by
$ ( y^2+zx : =  x_0,   x^2 :  = x_1, yz : = x_2, xy : = x_3) $
( in the projective embedding given by $ (zx, yz, x^2, xy, y^2)$
it corresponds to projection from the point $(1,0,0,0,-1) \in \PP^4
\setminus Y$).

In this case one sees that the matrix $\cal B$ is
$$\cal B=\begin{pmatrix}
x_3 &-x_2 &-x_0\\ -x_1 &0 & x_3\\ 0 &x_3 &-x_1
\end{pmatrix},$$
     and that $G$ is then the cubic of equation
$-x_3^3+x_1^2x_2+x_1x_2x_3=0$.

\end{rem}

\begin{defn}
We define now the {\bf  direct construction } of the bundle $\E$
relying on our results above.

Consider a sheaf $\cal G$ defined by an exact sequence as in \ref{Geq1},
and which is invertible on a cubic surface $G$ (i.e., at each point
$ y \in \PP^3$ $ rank ( \cal G \otimes \C_y ) \leq 1$).

Define  $ \tau : = \cal G^{\otimes2}(-1)$ and let $\E$ be
a vector bundle which is an extension of $ 6 \hol$ by $\tau$
     as in \ref{seq1a} ( here and elsewhere, $\hol : = \hol_{\PP roj(V)}$ ).
\end{defn}

\begin{prop}\label{un2}
     $\E$ as above is unique up to isomorphism in the following cases:
\begin{enumerate}
\item
     if $G$ is a smooth cubic surface.

\item
      if $G$ is reducible to the union of a plane $T$ and
a smooth quadric $Q$ intersecting transversally.
\end{enumerate}
\end{prop}
\proof
As before, it suffices to show that
$\dim Ext^1(\tau,\cal O) = 6$.

Now, $Ext^1(\tau,\cal O)=H^0(\cal Ext^1(\tau,\cal O))$
and the exact sequence
$$
0\rar \cal O(-3) \rar \cal O \rar \cal O_G \rar 0
$$
yields
$$
0 \rar \cal Hom(\tau, \cal O_G) \rar \cal Ext^1(\tau, \cal O(-3))
\rar \cal Ext^1(\tau, \cal O)
$$
where the last map is $0$ on $G$.
Hence, $\cal Ext^1(\tau,\cal O)\cong \cal Hom(\tau, \cal
O_G(3))=\tau^\vee(3)$.

Using the previous notation for the Cartier divisors corresponding to
$ \cal G$ and $ \hol_G(1)$, we want to show that
$ h^0 ( \hol_G (4 H - 2 L ) = 6 $.

Assume first that $G$ is a smooth cubic surface:
then by Riemann Roch it suffices to show the vanishing
of the first cohomology group $ h^1 ( \hol_G (4 H - 2 L ) = 0 $.

We argue as before using  Ramanujam's vanishing theorem,
since $$ h^1 ( \hol_G (4 H - 2 L ) =
     h^1 ( \hol_G (- 5 H + 2 L ) = h^1 ( \hol_G (- 13 L + 5 \sum_i E_i )),$$
and $ |13 L - 5 \sum_i E_i | \supset |10 L - 4 \sum_i E_i | + |H|$
contains a reduced and connected divisor.

In the second case, observe that there is no birational morphism
of a smooth quadric $Q$ onto the plane, thus $\cal G$
defines $ \pi $ which is an isomorphism on the plane, and has degree
zero on  $Q$.

Since we know that $\pi$ is an embedding of $Q \cap T $,
$\pi |_Q$ is the projection of $ Q \cong \PP^1 \times \PP^1$
on the second factor, followed by the isomorphism
of $\PP^1 $ onto $Q \cap T $.

Since $\tau^\vee(3)\mid_T=\cal O_T(2)$,
$\tau^\vee(3)\mid_Q=\cal O_Q(4,0)$,
and  $ H^0 (\tau^\vee(3)\mid_T=\cal O_T(2)) \ra
H^0 (\tau^\vee(3)\mid_Q=\cal O_Q(4,0))$ is surjective
we obtain
$$
\dim Ext^1(\tau,\cal O)=h^0(\tau^\vee(3))=6.
$$

\QED

\begin{rem}
Indeed the above proof shows that if $ G = T \cup Q$,
with $Q$ a smooth quadric, and $\cal G$ is invertible,
then necessarily $T$ and $Q$ intersect transversally.
\end{rem}

We observe now that Lemma \ref{sym} provides a resolution of
$ \tau : = \cal G^{\otimes2}(-1)$  starting
from \ref{Geq1}.
We take the second symmetric power of the sequence (\ref{Geq1})
and we obtain a resolution:
\begin{equation}
         \label{eq:T}
         0\rar (\Lambda^2 W') \otimes \cal O(-2) \rxar{\neg \cal B(-1)}
(U'\otimes W') \otimes \cal O(-1)
\rxar{\tilde{\cal B}} (S^2 U')\otimes \cal O
\rar
\tau(1)\cong\cal G^{\otimes 2} \rar 0
\end{equation}
where ${\neg \cal B}$ is the contraction given by the composition
of the natural inclusion from
$(\Lambda^2 W') \otimes \hol(-1)$ to $(W'\otimes W')  \otimes \hol(-1)$
with the map $\cal B\otimes id_{W'}(-1)$,
while $\tilde {\cal B}$ is the composition of $ id_{U'} \otimes \cal B$
with the surjection  $(U' \otimes U') \otimes \hol$ to $(S^2 U')
\otimes \hol$.

Consider now the exact sequence defining $\E$
$$ 0 \rar 6 \cal O \rar \cal E \rar \tau \rar 0,$$
and the above projective resolution of $\tau$:
by the mapping cone construction (cf. e.g. \cite{ei}, pages 650-651)
we obtain a projective resolution of
$\cal E$:
\begin{equation}
         \label{eq:resE}
         0\rar \Lambda^2 W' \otimes\cal O(-3) \rxar {\neg\cal B(-2) }
U'\otimes W'  \otimes \cal O(-2)
\rxar{(\tilde{\cal B}(-1),\lambda)}
\begin{matrix}
      6\cal O \\ \oplus \\ (S^2 U'  \otimes \cal O(-1))
\end{matrix}
\rar \cal E \rar 0.
\end{equation}

We now want to find a relation between the multiplication
map $B:\ U\otimes V \rar W$,
where $U$ (resp. $W$) denotes as usual $H^1(\cal E(-2))$ (resp.
$H^1(\cal E(-1))$),
and the above map $\tilde B:\ W' \rar U'\otimes V$.
Let $\cal E$ be the unique sheaf given by $B$ (cf. page 11).

We split the above resolution \ref{eq:resE} of $\cal E$ into two
short exact sequences,
denoting by $\cal K$ the image of $(\tilde{\cal B}(-1),\lambda)$.
This gives
$H^1_*(\cal E)\cong H^2_*(\cal K)=
\ker [H^3_*(\Lambda^2 W' \otimes\cal O(-3))
\rar H^3_*(U'\otimes W' \otimes \cal O(-2))]$.
Fixed these isomorphisms,
we can perform the following identifications:
$U\cong\ker[\neg\cal B:\ \Lambda^2 W'\otimes V^\vee \rar U'\otimes W']$,
$W=\Lambda^2 W'$,
and the multiplication map $B$ is given as the composition,
     in the following diagram
$$\xymatrix{
U\otimes V \ar@{^{(}->}[d] \ar[r]^B & W=\Lambda^2W'\\
\Lambda^2W'\otimes V^\vee\otimes V \ar[ur]_{\neg}},
$$
of the natural  inclusion with the natural contraction $\neg:\
V^\vee\otimes V \rar
\C$.

One can obtain the above factorization also in the
following alternative way :
Beilinson's complex yields the following short exact sequence:
$$0 \rar U\otimes \Omega^2(2) \rar W\otimes
\Omega^1(1)\oplus 6\cal O \rar \cal E \rar 0,$$
where $U$ (resp. $W$) denotes as usual $H^1(\cal E(-2))$ (resp.
$H^1(\cal E(-1))$).

    From the above we get:
\begin{equation}
      \label{eq:resE1}
0\rar
\begin{matrix}
      U\otimes \hol(-2)\\\oplus \\ W\otimes \hol(-3)
\end{matrix} \rar
\begin{matrix}
      (U\otimes V\otimes \hol(-1))\\
\oplus\\
(W\otimes \Lambda^3V\otimes \hol(-2))
\end{matrix} \rar
\begin{matrix}
      (6\hol)\\
\oplus\\
(W\otimes\Lambda^2V\otimes\hol(-1))
\end{matrix} \rar
\cal E \rar 0.
\end{equation}

Comparing the two sequences, we obtain the following identifications:
\begin{itemize}
\item $W'\cong \Lambda^2 W$,
\item $U\cong \ker[{\neg \cal B}:\ \Lambda^2 W'\otimes V^\vee \rar
U'\otimes W']$,
\item $U'\otimes W' \cong ( W\otimes \Lambda^3 V)/U$,
\item $S^2U'\cong (W\otimes \Lambda^2V)/(U\otimes V)$.
\end{itemize}

Based on the above considerations we give the following
\begin{defprop}\label{cross}

\ \\
The {\bf Cross-Product-Involution on Tensors of type 3 x 3 x 4}
is given as follows:

to a 5-uple $(U',W',V,\delta ',\cal B)$,
where
$U', W'$ are 3-dimensional vector spaces, $V$ is a 4-dimensional vector space,
$\cal B\in Hom(W'\otimes V^\vee,U')={W'}^\vee\otimes V\otimes U'$, $\delta ':\
\Lambda^3W'\cong \CZ$ an isomorphism,
     we associate the 5-uple
$(U,W,V^\vee,\delta ,B)$, where:
\begin{enumerate}
\item $W:=\Lambda^2W'$ and,
since then $W$ is canonically isomorphic to ${W'}^\vee$, by
the duality $W'\otimes \Lambda^2 W' \rar \C $ induced by $\delta '$,
we let $ \delta : = {\delta '}^\vee$,

\item $U:=\ker[{\tilde B}:\ \Lambda^2 W'\otimes V^\vee \rar U'\otimes W']$,
where ${\tilde B}$ is defined as above through the contraction ${\neg \cal B}$;
\item $B\in Hom(U\otimes V,W)=W\otimes V^\vee\otimes U^\vee$,
is the composition of the  inclusion $U\otimes V\hookrightarrow
\Lambda^2 W'\otimes V^\vee\otimes V$ with the natural contraction $\neg$.
\end{enumerate}
The dimension of $U$ is equal to three if we make  the

{\bf MAIN ASSUMPTION: }
${\tilde B}$ is surjective (this in turn  obviously implies the
injectivity of the map
$\cal B:\ {U'}^\vee \rar V\otimes W$).

The cross-product involution is then defined through the
associated 5-uple on the open set of tensors
satisfying the main assumption, and it is involutive
whenever the composition is defined.
\end{defprop}

\begin{proof}

Given a  5-uple $(U',W',V,\delta,\cal B)$, let
$(U,W,V^\vee,\delta,B)$ be  its corresponding 5-uple, to which
corresponds a third 5-uple
$(U'',W'',V,\delta,B'')$.
We have: $\cal B\in {W'}^\vee\otimes V\otimes U'$,
$B\in W \otimes V^\vee \otimes U^\vee$,
$B''\in W\otimes V\otimes {(U'')}^\vee$.

We claim that there exists a canonical isomorphism $U'={(U'')}^\vee$,
equivalently, a canonical isomorphism $U''={(U')}^\vee$.

     To show this, we shall first observe that both spaces can
canonically be regarded as subspaces of $ W \otimes V$,
and then, since both spaces have
the same dimension (for $U''$, this is a consequence of the hypothesis that
$(U,W,V^\vee,\delta,B)$ also satisfies the main assumption), it will suffice
to show that $ {(U')}^\vee \subset  U ''$.

$U''$ is the kernel of $\neg B:\ \Lambda^2 W'\otimes V \rar U^\vee\otimes W'$.
Identifying $\Lambda^2 W'$ with $W$,
the previous map becomes $\neg B:\ W\otimes V \rar U^\vee\otimes W^\vee$.
We consider now $\cal B$ as the map $\cal B:\ {U'}^\vee \rar V\otimes W$.
It suffices now to show $(\neg B)\circ \cal B ({U'}^\vee)=0$, i.e.,
by dualizing,
that $\cal B^\vee \circ (\neg B)^\vee\, (U\otimes W)=0$.
This is a consequence of the commutativity of the following diagram:

$$\xymatrix{
U\otimes W  \ar@{^{(}->}[d] \ar[r]^{\neg B}
& V^\vee\otimes W^\vee \ar[r]^{\cal B}
& U'\\
(\Lambda^2W'\otimes V^\vee)\otimes W
     \ar[ur]_\neg \ar[r]_{\neg \cal B\otimes id_W}
&(U'\otimes W')\otimes W   \ar[ur]_{\neg}
}.$$
\end{proof}

\begin{rem}
The two tensors considered in remark \ref{doubleline},
whose respective determinants yield the two non normal irreducible cubics
(which are not projectively
equivalent) satisfy
the main assumption. But the crossproduct involution constructs out
of them two tensors which do no longer satisfy the main assumption,
and which are projectively equivalent:

$$\begin{pmatrix}
x_2^\vee &x_3^\vee &0\\
0 &x_2^\vee &x_3^\vee\\
0&0&0
\end{pmatrix} \qquad
\begin{pmatrix}
0&0&0\\
x_0^\vee &x_2^\vee &0\\
x_2^\vee &0 &x_3^\vee
\end{pmatrix}.$$

\end{rem}

\section{Semistability of $\E$ and Moduli}

In this section we try to have a broader outlook at the
vector bundles
$\cal E$ that we are considering in this paper.
We shall show that their explicit geometric construction lends itself
to construct a
natural moduli space $\frak A ^0$ for them. Since moduli
space for vector
bundles have been constructed in great generality by Maruyama, it
seems natural to
investigate their Gieseker stability (we refer  to \cite{oss} and
especially to \cite{hl}
as  general references). We
conjecture that our bundles are   Gieseker stable, but  unfortunately
for the time being
we only managed to prove  their  slope (Mumford-Takemoto) semistability.

We are however able to prove that our vector bundles are simple,
and we observe then (cf. Theorem
2.1 of \cite{kob}) that moduli spaces of simple vector bundles exist 
as (possibly
non Hausdorff) complex analytic spaces.

We show indeed that the above moduli space exists as an algebraic variety.
  More precisely
that,  under a suitable open condition, we can construct a G.I.T. 
quotient $\frak A
^0$ which is a coarse moduli space.

Recall that , by \ref{ChernC}
 a vector bundle  $\E$  be given as in \ref{E=ker}
or as in \ref{E=ker2}
has  total Chern class 
\begin{equation}\label{chern(t)}
     c(\E) (t) = 1 + 3 t + 6 t^2 + 4 t^3.
\end{equation}

The next lemma will lead to a characterization of the vector bundles
obtained from the
direct construction as an open set in any family of
vector bundles with the above Chern polynomial.

\begin{lem}\label{h0=6}
Consider a  rank 6 vector bundle of  $\E$ with
     total Chern class
\begin{equation*}
     c(\E) (t) = 1 + 3 t + 6 t^2 + 4 t^3.
\end{equation*}
and assume that $ H^i(\E) = 0 \ \forall i \geq 1$:
then $h^0(\E) = 6$.

\end{lem}

\proof
Under these assumptions  $h^0(\E) = \chi (\E) $,
$\chi(E)$ is determined from the Riemann Roch theorem, and
we know that there are vector bundles $\E$ for which
$ H^i(\E) = 0 \ \forall i \geq 1$ and
     $h^0(\E) = 6$ (cf. lemma \ref{expected}).

\qed

\begin{prop}\label{sstable}
Consider a  rank 6 vector bundle of  $\E$ with
     total Chern class
$1 + 3 t + 6 t^2 + 4 t^3$,
such that
\begin{enumerate}
\item
$h^0(\E) = 6$
\item
the 6 sections generate a rank 6 trivial subsheaf
with quotient $\tau$
\item
$h^0(\E^\vee) = 0$
\item
$ \E$ is a subbundle of $ 3 \Omega^1 (2)$.
\end{enumerate}
Then $\E$ is  slope-semistable.
     \end{prop}

\proof

Let $\E''$ be a destabilizing subsheaf of rank $r \leq 5$ and
maximal slope $\mu = d/r$ :
without loss of generality we may assume that $\E''$ is  is a
saturated reflexive subsheaf,
and similarly
$\tilde{\E} : = \E'' \cap 6 \hol$ is a saturated reflexive subsheaf
of $6 \hol$.
$$
\xymatrix{
0 \ar[r] &6\hol \ar[r] &\E \ar[r] &\tau \ar[r] &0\\
0 \ar[r] &\tilde\E \ar[r] \ar[u] &\E'' \ar[r] \ar[u] &\tau'' \ar[r] \ar[u] &0\\
&0 \ar[u] &0 \ar[u]\\
}
$$

The slope $ \mu(\E) $ equals $ 1/2$. On the other hand, by hypothesis 4
and since  $\Omega^1 (2)$ is a stable bundle (cf. 1.2.6 b , page 167 of
\cite{oss}), the slope of $\E''$ is at most $2/3$, and
     $ < 2/3$ unless $\E'' \cong \Omega^1 (2)$.

{\bf CLAIM:} $\E$ contains no subsheaf isomorphic to $\Omega^1 (2)$.

{\em Proof of the claim:} $ h^0 (\Omega^1 (2)) = 6 = h^0 (\E)$,
thus this calculation contradicts hypothesis 2.

\qed

We have that $ d : = c_1 (\E'') = c_1 (\tilde{\E}) + c_1 (\tau'')$,
and $\tau'' \subset \tau$ is a coherent subsheaf supported on a divisor, thus,
     $ c_1 (\tau'')\leq   c_1 (\tau) = 3$.

On the other hand, $ c_1 (\tilde{\E}) \leq 0$, and if equality holds,
then  $ \tilde{\E} \cong r \hol$.

Hence, $ 1 \leq d  \leq 3$, and we have
$$  2/3 >  \mu = d/r > 1/2  \Leftrightarrow   4 d > 2 r > 3d . $$

These inequalities leave open only the case $ d=3, r=5$.

We show that this case does not exist.

In fact, otherwise we  consider the  quotient by the subbundle
$ \tilde{\E} \cong r \hol$. By hypothesis 3, and the proof of
lemma \ref{unicity} we see that $ \E / \tilde{\E}$ is
an extension corresponding to a tensor of maximal rank, hence
it yields a vector bundle $\cal V$ (cf. corollary \ref{vb}).

Since the torsion sheaf $\tau'' \subset \cal V$, we obtain
$\tau'' = 0$, hence $ d \leq 0$, absurd.

\QED

\begin{rem}
The possible exceptions to slope-stability, in view of the inequalities
$$  2/3 >  \mu = d/r  \geq 1/2  \Leftrightarrow   4 d  \geq 2 r > 3d  $$
are:

1. $  d=1, r=2 $

2. $  d=2, r=4 $.

Matei Toma pointed out how case 2. could be excluded using
Bogomolov's inequality for stable bundles, as done in Lemma 3.1 of his
paper \cite{tom}. The case $r=2, c_1 (\tilde{\E}) = -2$ seems 
as of now the most difficult case.

Observe that slope-stability of $\E$  implies Gieseker stability of
$\E$ , which in turn
implies that there is a point in the moduli space of Gieseker
semistable bundles
corresponding to the
isomorphism class  of $\E$.
\end{rem}

\begin{lem}\label{simple}
Let $\E$ be a vector bundle as in (\ref{E=ker})
with $h^0(\E)=6$ (equivalently, $h^1(\E)=0$)
    and verifying the second assumption.
Then $hom (\E, \E) = 1$, i.e., $\E$ is simple.
\end{lem}

\proof
We consider the exact sequence
$$0 \ra Hom(\E,6\hol) \ra Hom(\E,\E) \ra Hom(\E,\tau) \ra Ext^1(\E,\hol).$$

We have
$Ext^1(\E,\hol)\cong H^1(\E^\vee)\cong H^2(\E(-4))$
and from the exact sequence (\ref{E=ker}) we infer $ H^2(\E(-4)) = 0$.
Since
$Hom(\E,6\hol)=0$ by proposition 3.2, it follows that $Hom(\E,\E) \cong
Hom(\E,\tau)$.

We compute $hom(\E,\tau)$ by considering the exact sequence
$$0 \ra Hom(\tau,\tau) \ra Hom(\E,\tau) \ra Hom(6\hol,\tau).$$
Indeed $hom(\hol,\tau)=h^0(\tau)=0$ (since $h^0(\E)=6$) and, since $\tau$ is
$\hol_G$-invertible, we have $hom(\tau,\tau)=1$.

\qed

\begin{lem}\label{dim=19}
Let $\E$ be a  simple  vector bundle of rank 6, with Chern
classes $ c_1(\E) = 3, c_2(\E) = 6, c_3(\E) = 4$.
Then the local dimension of the  moduli space
$ \frak M^{s} ( 6; 3,6,4)$ of simple vector bundles
at the point corresponding to
$\E$  is at least $19$.
\end{lem}

\proof
The moduli space of simple vector bundles exists (cf. \cite{kob},
Theorem 2.1) and it is well known
that the local dimension is at least equal to the expected dimension $ h^1 (
\E^\vee
\otimes
\E) - h^2 (
\E^\vee
\otimes
\E)$. On the other hand, $\E$    simple means that $ h^0 ( \E^\vee \otimes \E)
= 1$, hence follows also that $ h^3 ( \E^\vee \otimes \E) = 0$,
since by Serre duality $ h^3 ( \E^\vee \otimes \E) =
     h^0 ( \E^\vee \otimes \E (-4)) = 0$.

Thus the
expected dimension  equals $ - \chi ( \E^\vee \otimes \E)  + 1 $ and there
     remains to calculate $ - \chi ( \E^\vee \otimes \E) $.
This can be easily calculated in the case where we have an exact  sequence
$ 0 \ra \E \ra 9 \hol (1) \ra 3 \hol(2)$. We omit the rest of the
easy calculation.

\qed

In the following theorem we shall phrase the geometric
meaning of the cross product involution in terms of a
birational duality of
moduli space of vector bundles, $\frak A ^0$ on $\PP^3$,
${\frak A ^0_*}$ on ${\PP^3}^\vee$.

{\bf Theorem C}
  {\em    Consider the moduli space $ \frak M^s ( 6; 3,6,4)$ of
  rank 6 simple vector
bundles $\cal E$ on $\Pn 3$
     with Chern polynomial $ 1 + 3t + 6 t^2 + 4 t^3$,
and inside it the  open set  $\ \frak A$  corresponding to the simple 
bundles with minimal cohomology, i.e., such that
\begin{enumerate}
\item
$ H^i(\E) = 0 \ \forall i \geq 1$
\item
$ H^i(\E (-1)) = 0 \ \forall i \neq 1$
\item
$ H^i(\E (-2)) = 0 \ \forall i \neq 1$.
\end{enumerate}

Then   $\frak A$  is irreducible of dimension 19
and it is bimeromorphic to $ \frak A^0 $,
where $\frak A^0$ is an open set of the G.I.T. quotient space of the
projective space $\frak B$ of
tensors of type $(3,4,3)$, $\frak B : = \{B\in \PP ({U}^\vee\otimes
V^\vee\otimes W)\}$
by the natural action of $SL(W) \times SL(U)$ (recall that $U,W$ are two
fixed vector spaces of dimension 3, while $ V = H^0( \PP^3, \hol (1))$.

Let  moreover $ [B] \in \frak A ^0$ be a general point:
then to $[B]$
corresponds a vector bundle $\E_B$ on $\PP^3$, and also a vector bundle
$\E^*_B$ on ${\PP^3}^\vee$, obtained from the direct construction
applied to $\cal G^*_B$ (cf. \ref{*}).  $\E^*_B$
is the vector bundle $\E_{\cal B}$, where $[{\cal B}] \in {\frak A ^0_*}$
is obtained from $B$ via the cross product involution.}

\proof
To any such tensor $B$ we tautologically associate two linear maps which
we denote by the same symbol,
$$B :  U \otimes V \ra W , \quad B :  U \otimes V \otimes  \hol(1)\ra W \otimes
\hol(1)$$
and using the Euler sequence we define a coherent sheaf $\E$ on $\PP^3$
as a kernel, exactly as in the exact sequence \ref{E=ker2} (except that
surjectivity holds only for $B$ general).

As we already saw in \ref{B=beta}, this is equivalent to giving
$\E$ as the kernel of a homomorphism $\beta$ as in \ref{E=ker}.
Observe that $GL(W) \times GL(U)$ acts on the vector space of such
tensors, preserving the isomorphism class of the sheaf thus obtained.

We define   $\frak B' $ as the open set in  $\frak B $ where
$\beta$ is surjective (thus $\E$ is a rank 6 bundle) and
$h^0(\E)=6$. Both conditions amount to the surjectivity
of $h^0(\beta)=h^0(B, \epsilon)$, cf. \ref{E=ker2}. 
We further define
$\frak B ''$ as the smaller open set where the second assumption 
is verified,
and we observe then that  lemma \ref{simple} 
ensures the existence of a
morphism
$\frak B '' \ra \frak A$ which factors through the action of
$ SL(W) \times SL(U)$.

Since we want to construct a G.I.T. quotient of an open set of $\frak B$,
we let $\frak B^*$ the open set of tensors $B$ whose determinant
defines a cubic surface $G^* \subset {\PP^3 }^\vee$, i.e., we have an
exact sequence on ${\PP^3 }^\vee$ of the form (set $\hol_* : =
\hol_{{\PP^3 }^\vee}$)
\begin{equation}\label{*}
         0 \rar U\otimes \cal O_*(-1)  \rxar{B} W\otimes \cal O_*  \rar
\cal G^* \rar 0.
\end{equation}
Since the determinant map is obviously $SL(W) \times SL(U)$-invariant,
the tensors in  $\frak B^*$ are automatically semistable points
for the $SL(W) \times SL(U)$-action, by virtue of the criterion of
Hilbert-Mumford.

Observe now that the maximal torus $\C^* \times \C^*$ of
$ GL(W) \times GL(U)$ acts trivially on $\frak B$, thus
we get an effective action of $SL(W) \times SL(U)$ only upon
dividing by a finite group $  K' \cong (\Z/3)^2$.

We claim that  $(SL(W) \times SL(U)) / K'$ acts freely on
the open subset $\frak B^{**}   \subset \frak B^* $,
$\frak B^{**} = \{ B \in\frak B^*  | End(\cal G^* ) = \C \}$.

This is clear since the stabilizer of $B$ corresponds uniquely
to the group of automorphisms of $ \cal G^* $,
and any such automorphism acts on $W \cong H^0(\cal G^* )$,
and induces a unique automorphism of $U$ in view of the exact
sequence \ref{*}.  But every automorphism is multiplication
by a constant, thus it yields an element in $K'$.

We want to show that the orbits are closed.
But the orbits  are contained in the
fibres of the determinant map: thus,
it suffices to show that, fixed the cubic surface $G^*$, if
we have a 1-parameter family where $ \cal G_t \cong \  \cal G_1$
for $ t \neq 0$, then also $ \cal G_0 \cong \  \cal G_1$.

This holds on the smaller open set $\frak B^{***} \subset \frak B^{**}$
consisting of the tensors such that the cubic surface $G^*$
is smooth: since then $ \cal G_0$ is invertible, and the Picard group
of  $G^*$ is discrete.

We have proven that $\frak B^{***}$ consists of stable points,
and observe that  the condition $End(\cal G^* ) = \C $
holds if $\cal G^* $ is $\hol_{G^*}$-invertible, or
it is torsion free and $ G^*$ is normal.
Therefore the open set $\frak B^{st}$ of stable points is nonempty.

We define $\frak A ^0$  as the open set of the G.I.T. quotient
corresponding to
$\frak B^{st} \cap \frak B ''$.

The fact that $\frak A$ is irreducible follows since
every bundle $\E$ in $\frak A$ has a cohomology table which
(by Beilinson's theorem, as in the beginning of section 3)
implies that $\E$ is obtained from a tensor
$B$ in the open subset ${\frak B ' }^0 \subset {\frak B ' }$
consisting of those $B$ for which the corresponding bundle $\E$ is simple
(note that ${\frak B ' }^0  \supset {\frak B '' }$).

Now,  $\dim \frak A^0 = 19 $, while $\dim \frak A \geq 19 $
by \ref{dim=19}; we only need to observe that if
$[B] , [B'] \in \frak A^0$ and two bundles
$\E_B$ and $\E_{B'}$ are isomorphic, then the corresponding
tensors  $B, B'$ are $ GL(U) \times GL(W)$ equivalent, since
they express the multiplication matrix for the intermediate cohomology
module $ H^1_*( \E)$. Thus $[B] = [B'] \in \frak A^0$.

It follows on the one side that
$\frak A^0 $  parametrizes isomorphism classes of bundles,
and on the other side  $\frak A^0 $
maps bijectively to an open set in $\frak A$, in particular
$\dim \frak A = 19 $, since $\frak A $ is irreducible.

\QED

\section{The explicit unirational family}\ \\

We have up to now studied extensively the vector bundles $\E$
such that an even set of $56$ nodes on a sextic surface
$F$ should come from a symmetric homorphism
associated to a section of $ S^2(\E)$.

We have however almost shown, because of  corollary \ref{21}
and of proposition \ref{gen} that all such sections have as determinant the
square of a cubic surface $G$, if the cubic $G$ appearing in
the direct construction is an irreducible cubic.

It seems therefore only natural to try to see what happens for a
reducible cubic,
hoping that then
$h^0(S^2\cal E) > 21$.

We assume henceforth that $G$ is the union of a smooth quadric $Q$
with a plane $T$ intersecting transversally.
We have already observed in the proof of proposition
    \ref{un2} that in this case there is
a unique choice for $\cal G$ , likewise for $\E$.

\begin{lem}\label{22red}
If $G$ is the union of a smooth quadric $Q$
and of  a plane $T$  which intersect transversally,  then $h^0
(\tau^{\otimes2}) = 1$ and $h^0(S^2\E)=22$.
\end{lem}

\proof
In this case the sheaf $\cal G$
corresponds to the sheaf $\hol_Q(0,2)$ on $Q$ and to
$\hol_T(1)$ on $T$.

Therefore $\tau:=\cal G^{\otimes 2}(-1)$ corresponds to
$\hol_Q(-1,3)$ on $Q$ and to $\hol_P(1)$ on $P$.

    Thus the sections
of $H^0 (\tau^{\otimes2})$ vanish identically on $Q$ and correspond
to section of $\hol_T(2)$  vanishing on $ Q \cap T$.

The second statement follows then from the proof of corollary \ref{21}.

\qed

We also remark that, since we assume $Q\cap T$ is smooth,  $G$ is
unique up  to projective equivalence.

\bigskip

We shall now give an explicit tensor $B_0$ whose associated sheaf
is the unique $\cal G$ on a reducible cubic of the form $ T \cup Q$,
where $Q$ and $Q\cap T$ are smooth, and compute explicitly that the tensor
corresponding to the
unique $\E$ gotten from the direct construction is again $B$:
this will allow us to calculate explicitly the determinant of a generic
symmetric map
$ \E^\vee \ra \E$, and to show that it is a nodal sextic.

\begin{lem}
Consider the following $3\times 3\times4$ tensor $B_0$:
\begin{equation}\label{normalForm}
B_0 =
\begin{pmatrix}
0 & -x_3 & x_2 \\
x_3 & 0 & -x_1 \\
-x_2 & x_1 & 0\\
\end{pmatrix}
+ x_0
\begin{pmatrix}
1 & 0 & 0  \\
0 & 1 & 0 \\
0 & 0 & 1\\
\end{pmatrix}.
\end{equation}
The sheaf $\cal G_0$ associated to $B_0$ is an invertible
sheaf on the reducible cubic $G_0: = \{ x_0 ( \sum_i x_i^2 ) = 0 \} $.

Its class in $\frak A ^0$ is invariant under the cross product involution.
More precisely, the vector bundle $\E_0$  obtained from $B_0$
via the direct construction
   has the required cohomology table and it has again $B_0$
as multiplication matrix for the intermediate cohomology
module $M_0 : = H^1_*(\E_0)$.

Moreover $h^0(S^2\E_0)=22$.
\end{lem}

\proof
The determinant of $B_0$ equals $G_0$.
On the plane $x_0=0$ the Pfaffians
are $x_1=x_2=x_3=0$ and $\mu$ has rank 2.
Elsewhere $G_0$ is smooth, whence $\cal G$ is
everywhere invertible.

We show now that its class is invariant under the cross product
involution.
Denote by $W'$ and $U'$ the vector spaces of rows and colums of $B_0$.

Then $U$ is defined as the kernel of the composition
${\neg B_0}:\ \Lambda^2 W'\otimes V^\vee \rar (W'\otimes V^\vee)\otimes W'
\ra U'\otimes W'$.
Let $e_1,e_2,e_3$ (resp. $f_1,f_2,f_3$) be the standard basis of
$W'\cong\C^3$ (resp. $U'$).
Then we choose $e_2\wedge e_3, -e_1\wedge e_3,e_1\wedge e_2$
as basis of $\Lambda^2 W'\subset W'\otimes W'$
and  $f_1\otimes e_1,f_2\otimes e_1,\ldots,e_3\otimes f_3$
as basis of $U'\otimes W'$ and, viewing the map ${\neg B_0}$
as a $ 9 \times 3$ matrix of linear forms on the space  $V^\vee$, we have:
$$\tilde{B_0}=\begin{pmatrix}
0 &(B_0)^3 &-(B_0)^2\\
-(B_0)^3 &0 &(B_0)^1\\
(B_0)^2 &-(B_0)^1 &0
\end{pmatrix},$$
where $(B_0)^i$ denotes the $i$-th column of $B_0$.
The kernel $U$ is then:
$$U=\{ u=\sum a_i x_i^\vee| a_i\in \C^3,\  \tilde{B_0} \, u=0\}=
\langle
\begin{pmatrix}
x_0^\vee\\ x_3^\vee \\ -x_2^\vee
\end{pmatrix},
\begin{pmatrix}
-x_3^\vee \\ x_0^\vee \\ x_1^\vee
\end{pmatrix},
\begin{pmatrix}
x_2^\vee \\ -x_1^\vee \\ x_0^\vee
\end{pmatrix}
\rangle.$$

A rapid inspection shows that, if we choose the above bases and the dual
basis of $V^\vee$, the $3\times 3 \times 4$ tensor that we obtain is exactly
identical to $B_0$.

The direct construction gives then a unique vector bundle $\cal E$
(cf. proposition \ref{un2}), and we claim that its  cohomology table  is the
required one. As in the proof  of lemma \ref{expected}
it suffices  to calculate the dimensions of the cohomology groups
$h^i(\tau (-n))$, for $ 0 \leq n \leq 3$ and to show their vanishing for
$ i = 0, i = 2$.

We use the exact sequence
   $$0 \ra \tau \ra \hol_Q(-1,3) \oplus \hol_T(1) \ra
\hol_{Q \cap T}(1) \cong \hol_{\PP^1}(2) \ra
0$$
which is soon seen to be exact on global sections, hence
the case $i=0$ follows right away.
For the case $ i=2$ we observe that in the exact cohomology sequence
   $$0 \ra  H^1(\tau(-n)) \ra H^1 ( \hol_Q(-1-n,3-n) ) \rxar \psi
H^1 (\hol_{Q \cap T}(1-n) ) \ra H^2(\tau(-n)) \ra 0
$$
$\psi$ is surjective, since its cokernel
is isomorphic to $ H^2 ( \hol_Q(-2-n,2-n) ) $
whose dimension equals $ h^0 ( \hol_Q(n, n-4 ) ) = 0$
(since $n < 4$).

The last claim follows from lemma \ref{22red}.

\QED

\begin{rem}

Consider the invertible sheaf given by:
\begin{equation}\label{MultTab}
0\rar 3\cal O(-1) \rxar {B_0} 3\cal O \rar \cal G_0\rar 0.
\end{equation}

We have that $h^0(\cal G_0)=3$ and the associated morphism  $\pi$ from
$G_0$ to
$\Pn 2$ is determined by the rational map given by
the entries of any column of $Ad(B_0)$:
$$
Ad(B_0)=
\begin{pmatrix}
         x_1^2+x_0^2 & &\\
         x_1x_2+x_0x_3 &x_2^2+x_0^2 &\\
         x_1x_3+x_0x_2 &x_2x_3+x_0x_1 &
x_3^2+x_0^2
\end{pmatrix}.
$$

We see right away that $\pi$ is the identity on
the plane $T$,
and the projection along one ruling from the quadric $Q$ to $Q\cap T$.

\end{rem}

\begin{prop}\label{exist}
For a section  $\phi \in H^0(S^2\E_0)$, denote by
$\varphi\in Hom(\E^\vee,\E)$ its associated symmetric morphism.
Then, for $\phi$ general, $F:=\{x\mid\det(\varphi)=0\}$ is a nodal sextic
surface with, as singularities, exactly an even set of 56 nodes
$\Delta=\{x\mid Corank(\varphi) =  2\}$.
\end{prop}

\proof

  The required computations were performed  and can be performed
and verified by using the computer-algebra system
\cite{M2} over a finite field, or  over $\QQ$
(cf.  script I in the Appendix).

The first step is to compute explicitly the fibre over $B_0$ inside
the variety of pairs $\frak M_{AB}$ (cf. \ref{AB}),
i.e., the vector space of symmetric matrices
$A\in Mat(12\times 12,H^0(\hol_\Pn 3(2)))$
satisfying the equation $(B_0,\epsilon)\cdot A=0$.

Step two: for a random $A$ in such a fibre one computes  the
G.C.D. of  two (different) $6\times 6$ minors of $A$:
if the G.C.D.  has degree 6, then it is the equation of the sextic $F$.

Step three: one verifies with the jacobian criterion that the singular
locus of
$F$  consists exactly of  a 0-dimensional subscheme of  length $56$.

Step four: one verifies that the ideal sheaf of the singular locus is
a radical ideal. Then the singularities are just a set of nodes.

A further (but not absolutely necessary) check consists in  verifying
that  the scheme $\Delta$
  coincides with  the subscheme formed by those $56$ reduced points:
this can be performed by computing a set of $5\times 5$ minors
of $A$ sufficient in order to generate the ideal of the 56 points.

\qed

\begin{rem}
Since the space of reducible cubic surfaces has dimension 12 (9 + 3),
we obtain an explicit  family parametrized by a rational variety
$\Phi_0$  of
dimension 33= 21 + 12.
\end{rem}
\proof
We simply construct a  parameter space by choosing
a 12-dimensional subgroup  $ H \subset \PP GL(V)$ such that
the orbit of $G_0$ dominates the space of reducible cubics,
and then we take as parameter space $ H \times \PP(H^0(S^2(\E)))$.

Then to the pair $ (g, \phi_0)$ corresponds the vector bundle
$g^*(\E ) : = \E_{g G_0}$ and the section $ g^*(\phi_0 )$,
and, correspondingly, the sextic surface $ g^*( det (\varphi)=0)$.

\qed

\begin{lem}
The morphism $ \Phi_0 \ra \PP (S^6(V)) $ associating
to $(g , \varphi)  $  the corresponding nodal sextic
$ F = det (\varphi)$ has fibres of dimension 6.
\end{lem}

\proof
Recall first that we have already shown that the surface $ g (G_0)$
uniquely determines a vector bundle $\E$ and conversely.

Second, observe that if $F$ has exactly $56$ nodes then $F$
determines the quadratic sheaf
$\F$ uniquely (observe  moreover that $\F$ has only constant
automorphisms).

Suppose   that  there are two different  vector bundles $\E, \E'$
and respective morphisms
$\varphi,\varphi'$ forming exact sequences as in theorem
\ref{CaCaThm} which define isomorphic cokernels  $\F , \F'$.

By abuse of notation we identify $\F '$ with $ \F$
and assume that we have $ \gamma : \E \ra \F$,
$ \gamma' : \E' \ra \F$ inducing such isomorphisms
of  the respective cokernels with $\F$.

A first question is whether there do exist lifts $\alpha:\ \E' \ra \E$ and
$\beta:\ {\E'}^\vee\ra \E^\vee$ of  the identity $id_{\F}$ on
$\F$ such that the following diagram commutes

\begin{equation}\label{injdiag}
\xymatrix{
0 \ar[r] &\E^{\vee} \ar[r]^\varphi &\E \ar[r]^\gamma &\F \ar[r] &0\\
0 \ar[r] &
{\E'}^{\vee} \ar[r]^{\varphi'} \ar[u]_\beta &
\E' \ar[r]^{\gamma'}\ar[u]_\alpha
&\F \ar[r]\ar[u]_{id_{\F}} &0\\ }.
\end{equation}

If such an $\alpha$ exists, then necessarily
the  submodule
$M : =  Im (H^1_*(\gamma))$ of $H^1_*(\F)$ equals
the  submodule
$M' : =  Im (H^1_*(\gamma'))$.

Assume now that $ M'= M$: then  since any automorphism
of $M$ lifts to an isomorphism of two minimal resolutions
of $M$, this automorphism induces an isomorphism
$\alpha$ of the respective first syzygy bundles, here
$\E$, resp. $\E'$ (cf.  the first page  of section 2).
We can henceforth assume that if $ M = M'$, then
$\E = \E'$, and then, using
   $Hom(\E,\E)=\C$ (cf. lemma \ref{simple}) we
get that $\alpha$ is the multiplication by a constant,
necessarily $\neq 0$.

  From the exact sequence
$$0 \ra Hom(\E^\vee,\E^\vee) \ra Hom(\E^\vee,\E)
\rxar{\gamma\circ} Hom(\E^\vee,\F),$$
it follows that there exists a unique homomorphism $\beta$
making the diagram commute. Again, using $Hom(\E,\E)=\CZ$,
we get that $\beta$ is the multiplication by a non-zero constant
and we have thus shown that if $M = M'$ then the sections
$ \phi$ and $ \phi '$ are proportional.

On the other hand, the choice of $M$ is completely
determined by the choice of a Lagrangian subspace $U$ of
the 6-dimensional space $ H^1 (\F(1))$, and we saw that for each
choice of $U$ there is a bundle $\E$ and a $\phi$
yielding an  exact sequence as in  theorem \ref{CaCaThm},
with  $M $ equal to the image of $ (H^1_*(\E))$.

We are done, since the dimension of the Lagrangian
Grassmannian $ LGr(3,6)$ equals 6.

\QED

We want to show that the explicit unirational family that
we constructed is locally maximal. To this purpose,
observe that to a pair $ (g, \phi_0)$ corresponds
a vector bundle $g^*(\E_0)$ and a section $g^*(\phi_0)$,
but more precisely a tensor $g*(B)$ and a matrix
of quadratic forms $A_{g, \phi_0}$ representing
$g^*(\phi_0)$ as in \ref{E=ker2}.

Thus $\Phi_0$ maps in a generically finite way to the 
variety of pairs  $\frak M_{AB}$ (cf. \ref{AB})
and we can consider the 
$GL(U) \times GL(W)$-orbit  of
its image.

Observe then that we obtain an irreducible algebraic set
$\Psi_0$ of dimension
$33 +1 +  9 + 9 - 1 = 51$.

The following lemma shows that $\Psi_0$ is indeed a component
of $\frak M_{AB}$.

\begin{lem}\label{TanAB}
Let $(B_0,A_0)\in\frak M_{AB}$ be a general point of the fibre over $B_0$.
Then the tangent space to $\frak M_{AB}$ at the point $(B_0,A_0)$
has dimension 51.
\end{lem}

\proof
Fixed the pair $(B_0,A_0)\in\frak M_{AB}$, we search the solutions
for a generic pair
$(B,A)\in Mat(3,12,\C)
\times
Mat_{Sym}(12,12,H^0(\hol_{\PP^3}(2)))$
of the equations
$$(B_0+tB,\epsilon)(A_0+tA) \equiv 0 \quad \mod{t^2}.$$

The above equation is equivalent to the two equations $B A_0 +B_0 A=0$,
  $\epsilon A=0$,
and we  have to compute the space of solutions.

Again we perform the computation by means of computer-algebra over
a finite field, it suffices to choose a point $A_0$ at random for which
the tangent space has dimension 51. The computation
works out  successfully (cf. script II in the Appendix).

\qed

We can now summarize the result of the construction of the above
explicit family:

\medskip

{\bf Main Theorem B.}
{\em There is a  family of nodal sextic surfaces
  with $56$ nodes forming an even set, parametrized by  a
smooth irreducible rational variety $\Phi_0$ of dimension
33, whose image $\Xi_0$ is a unirational subvariety  of dimension 27
of the space of sextic surfaces.
Moreover, the above family is versal, thus $\Xi_0$ yields
an irreducible  component of the subvariety of  nodal sextic surfaces with 56
nodes.}

\proof
The first assertions were proven in the sequel between lemma 6.1
and lemma 6.6

Let $\Xi$ be the subvariety of  nodal sextic surfaces with 56
nodes: since the property that the set of nodes is even is 
a topological property (cf.  for instance \cite{Ca}, \cite{Ca2}), 
it follows that there is an open and closed set
$\Xi' \subset \Xi$ such that for $ F \in \Xi'$ the set of 56 nodes
is even. We need only to prove that $\Xi_0 \subset \Xi'$ is open.

But  $\Xi' $  contains the open set $\Xi''$ such that, for
$F \in \Xi''$,  $H^1 (\F(2))$ has dimension 3 and the first assumption
is verified.

We can form a variety
$\Psi'$ consisting  of quadruples $ (F, U, B, \phi )$ where:

i)  $F \in \Xi' $ is a sextic surface,

ii) $U \subset H^1 (\F(1)) $ is a Lagrangian subspace, and,
$M$ being the intermediate cohomology submodule of
$ H^1_* (\F)$ determined by the choice of $U$ as in \ref{M},

iii) $B$ is a multiplication tensor for $M$ (depending on two choice
of bases)

iv) if $\E$ is the unique vector bundle determined
by $B$ as in \ref{E=ker2},
$\phi$
is a section of the vector bundle $ S^2 (\E)$ such that
$ det ( \varphi) = F$.

Then we see that the map $ \Psi' \ra \frak M_{A,B}$
is an embedding. Now lemma \ref{TanAB} shows that
$ \Psi_0 \subset \Psi'$ is open, and we are done.

\qed

\bigskip

It is a natural question to ask if the above is the unique
irreducible component of
the subvariety of  nodal sextic surfaces with  56
nodes forming an even set. 
For this purpose one should first settle the
case of Hilbert function (3,4) for $M$.

\section{The random approach}\ \\ \label {randomTrials}

Let $\mathbb M$ be a variety defined over a finite field of order $q$ and
let $\mathbb M_0\subset \mathbb M$ be a subvariety of codimension $k$.
The random approach consists in finding a point in $\mathbb M_0$ by
choosing at random points in $\mathbb M$.
Since the probability of hitting a point of $\mathbb M_0$ is $q^{-k}$
it is evident that this method is only successful if the computational time to
decide wether a point of
$\mathbb M$  actually belongs to $\mathbb M_0$ is small enough
(cf. \cite{Sch}).

In this section we show how
this method was applied to find the first examples
of sextic surfaces with an even set of 56 nodes.

\medskip

Let $\A $ denote the coordinate ring of $\PP^3$ and let $B$
  be the multiplication matrix of the intermediate cohomology module $M$.
  If $B$ is
general, since $\E$ is a syzygy bundle
(cf. section 3), it follows  (cf.  \ref{b}, and \ref{eq:resE}) that
$M$ has a resolution of the form

\begin{equation}\label{33res}
0\lar  M \lar 3\A [2] \lar  9\A [1] \lar
6\A \oplus 6\A [-1] \lar 9\A [-2] \lar 3\A [-3] \lar 0.
\end{equation}

In an analogous way to the one followed after the exact sequence
\ref{E=ker2}
we get that the symmetric morphisms $\varphi: \cal E^* \rar \cal E$
are exactly induced by the symmetric morphisms $a :\ 9\cal O(-1)
\rar  9\cal O(1)$
such that $b\circ a =0$, according to the following diagram:
$$
\xymatrix{
0 & \ar[l] 3\cal O(2) &\ar[l]_{b} 9\cal O(1) \ar[l] & \cal E \ar[l]\\
0 \ar[r] & \ar[r]^{\!^t\!b} 3\cal O(-2) \ar[u]_{0} &\ar[r] 9\cal
O(-1) \ar[r] \ar[u]_{a} & \cal E^\lor \ar[u]_{\varphi}
}.
$$
It is clear that the replacement of $(A,B)$ with $(a,b)$ reduces
the memory required for computations.

Repeated random choices of $b$ allow
to find an $\E$ with:
$$h^0(S^2\cal E)=\dim\{a : 9\cal O(-1) \rar
9\cal O(1)\mid a =\!^t\! a,\ b\circ a=0\}\geq 22.$$

This property leads to the definition of $\mathbb M$ and $\mathbb M_0$.

\begin{defn}
Let $\mathbb M$ be the Zariski open set
\begin{equation*}
\begin{split}
\mathbb M:=\{&b:\ 9\cal O(1) \rar 3\cal O(2) \mid
M:=\coker b \text{ has a resolution as in \ref{33res} and}
\\ &\iota:\ 6\cal O \rar \cal E:=\cal Syz_1(M) \text{ is injective}\}
\end{split}
\end{equation*}
$$\mathbb M_0:=\{b\in \mathbb M \mid h^0(S^2\E)\geq 22\}$$
\end{defn}

We already remarked that  $\mathbb M$ is nonempty.
A resolution for $S^2\E$ is provided by the following lemma.

\begin{lem}\label{S2resol}
If $0 \rar A \rar B \rar C \rar \cal E \rar 0$
is an exact sequence of locally free sheaves,
then the following sequence is also exact:
$$
0 \rar S^2A \rar A\otimes B \rar \Lambda^2B\oplus (A\otimes C) \rar B\otimes C
\rar S^2 C \rar S^2 \cal E \rar 0.$$
\end{lem}

\begin{proof}
By hypothesis we have $0 \rar  (B/A) \rar C \rar \cal E\rar 0$.
Therefore we get $0 \rar \Lambda^2(B/A) \rar (B/A)\otimes C \rar S^2
C \rar S^2 \cal E \rar 0$.
Resolutions for $\Lambda^2(B/A)$ and $(B/A)\otimes C$ are standard, respectively
$ 0 \rar S^2 A \rar A\otimes B \rar \Lambda^2 B \rar \Lambda^2(B/A)\rar 0 $ and
$0 \rar A\otimes C \rar B\otimes C \rar (B/A)\otimes C \rar 0$.
The resolution for $S^2\cal E$ stated in the lemma is the mapping cone
of the previous resolutions.

\end{proof}

Hence it was guessed that the ``good''  locus has codimension 7:

\begin{prop}
The condition $h^0(S^2\cal E)\geq 22$ is
expected to hold on a codimension 7 algebraic subset of $\mathbb M$.
\end{prop}

\proof

By applying the previous Lemma to a minimal free resolution of $\cal E$
we get a (non necessarily minimal) free resolution of $S^2 \cal E$:
\begin{equation}\small
          0 \rar 6\cal O(-6) \rar 27\cal O(-5) \rxar{r_3}
\begin{matrix}
          18 \cal O(-4) \\
          \oplus\\
          18 \cal O(-3)\\
          \oplus\\
          36 \cal O(-4)
\end{matrix}
\rxar{r_2}
\begin{matrix}
          54 \cal O(-3)\\
          \oplus \\
          54 \cal O(-2)
\end{matrix}
\rxar{r_1}
\begin{matrix}
          21 \cal O(-2)\\
          \oplus \\
          36 \cal O(-1)\\
          \oplus \\
          21 \cal O
\end{matrix}
\rxar{r_0}
S^2 \cal E
\rar 0.
\end{equation}

  Denote by $\cal K_i$ the image of the map $r_i$,
split the above exact sequence into  short exact sequences
and look at the associated cohomology exact sequences.
From
$$0 \rar H^o(\cal K_1) \rar H^0(21\cal O) \rxar{S^2(\iota^0)}
H^0(S^2\cal E) \rar
H^1( \cal K_1) \rar 0:$$
  and since $S^2\iota^0$ is injective ($\iota:\
6\cal O\rar \cal E$ being injective) we get, using also the other
cohomology sequences,
$0=H^0(\cal K_1)\cong H^1(\cal K_2) \cong
H^2(\cal K_3)$ and
$H^0(S^2\cal E)/H^0(21\cal O)\cong H^1\cal K_1 \cong H^2 \cal K_2$.

  We have also the short exact sequence
$$H^2(\cal K_3) = 0  \rar H^3(6\cal O(-6))  \rar H^3(27 \cal O(-5)) \rar
H^3(\cal K_3) \rar 0,$$
hence $H^3(\cal K_3)$ has dimension $48$.

Finally, the  exact sequence:
$$0 \rar H^2(\cal K_2)  \rar H^3(\cal K_3) \rxar\alpha H^3(54\cal O(-4))
\rar H^3 (\cal K_2) \rar 0,$$
since $H^3(\cal K_2)\cong H^2(\cal K_1)\cong H^1(S^2 \cal E)$ yields:
\begin{equation}
          0 \rar \frac{H^0(S^2\cal E)}{H^0(21\cal O)} \rar \C^{48}
\rxar\alpha \C^{54}
\rar H^1(S^2\cal E) \rar 0.
\end{equation}

Therefore the condition that $H^0(S^2\cal E)\cong H^0(21\cal O)$
  is equivalent
to the condition that the linear map $\alpha$ has maximal rank,
and since  we know that this happens in general,
the condition $h^0(S^2\cal E)\geq 22$ holds
in a determinantal subscheme   of  $\mathbb M $
of expected codimension $54-48+1=7$.

\qed

Let $\frak M_{ab}$ be the variety, analogous to $\frak M_{AB}$
(cf. \ref{AB}), of pairs (b,a)  such that $ a = ^t a, ab = 0$.
Computations similar to the ones in lemma \ref{TanAB} verified
over a finite field that in
a random point $(b_0,a_0)$ the variety of pairs $\frak M_{ab}$ is smooth
of dimension 123.

A standard argument then ensures the existence of a lift
of the  pair $(b_0, a_0)$ from a finite field to a number field
(cf. \cite{Sch}).

\medskip

This random approach, and the remark that the space of
reducible cubic surfaces is a codimension 7 subvariety of the
projective space of cubic surfaces led then to the
explicit family constructed in the previous section.

\section{Appendix: Macaulay2 scripts}

{\footnotesize

\noindent {\bf Script I:}
\begin{verbatim}
kk=ZZ/101
R=kk[x_0..x_3]
Pa=kk[a_1..a_78]
RP=R**Pa

--Step 0: the matrix (B,\epsilon)
use R
B=matrix{{x_0,-x_3,x_2},{x_3,x_0,-x_1},{-x_2,x_1,x_0}}
vars R
betti (B=diff(vars R,B)**R^{-1})
betti (epsilon=(vars R)**id_(R^3))
betti (B1=B||epsilon)
betti res coker B1
betti res prune coker B1

--Step 1: computing the 22-dim morphisms sol
use RP
A=genericSymmetricMatrix(RP,a_1,12);
equation=flatten (substitute(B1,RP)*A);
equation=substitute(diff(substitute(vars Pa,RP),transpose equation),R);
betti equation
betti (sol=syz(equation,DegreeLimit=>1)) --22

--Step 2: taking a random morphism A1 and computing the surface
A1=flatten (sol*random(R^22,R^1));
A1=substitute(substitute(A,Pa),A1);
B1*A1==0 --verification that it is really in the fiber
m1=random(R^6,R^12)*A1*random(R^{12:-2},R^{6:-2});
m2=random(R^6,R^12)*A1*random(R^{12:-2},R^{6:-2});
F=gcd(det(m1),det(m2));
degree F --6

--Step 3:
decompose ideal F; --sestica irriducibile
SingF=ideal(jacobian ideal F)+ideal F;
degree SingF,codim SingF --(56, 3)

--Step 4:
SingFsat=saturate SingF; --equiv. to: SingFsat=saturate(SingF,ideal vars R);
SingFsat==radical SingFsat --true

--Further check:
n1=random(R^5,R^12)*A1*random(R^{12:-2},R^{5:-2});
I=ideal det n1
i=1
while (i<50) do (
       n1=random(R^5,R^12)*A1*random(R^{12:-2},R^{5:-2});
       I=I+ideal det n1;
       i=i+1;
       );
codim I, degree I --(3,56)
I=ideal mingens I;
I=saturate I;
I==SingFsat --true
\end{verbatim}

\noindent {\bf Script II:}
\begin{verbatim}
kk=ZZ/101
R=kk[x_0..x_3]
Pa=kk[a_1..a_78]
Pb=kk[b_1..b_36]
P=kk[b_1..b_36,a_1..a_78]
RP=R**P
use R
B=matrix{{x_0,-x_3,x_2},{x_3,x_0,-x_1},{-x_2,x_1,x_0}}
vars R
betti (B=diff(vars R,B)**R^{-1})
betti (epsilon=(vars R)**id_(R^3))
betti (B1=B||epsilon)
betti res coker presB1
betti res prune coker presB1
use RP
A=genericSymmetricMatrix(RP,a_1,12)
equation=flatten (substitute(B1,RP)*A);
equation=substitute(diff(substitute(vars Pa,RP),transpose equation),R);
betti equation
betti (sol=syz(equation,DegreeLimit=>1)) --22
A1=flatten (sol*random(R^22,R^1));
A1=substitute(substitute(A,Pa),A1);

use P
betti (Bgen=map(P^{3:1},P^{12:0},genericMatrix(P,b_1,3,12)))
betti (B2=substitute(Bgen,RP)||substitute(epsilon,RP))
A2=A
equation=substitute(B1,RP)*A2-B2*substitute(A1,RP);
eq1=flatten equation^{0..2};eq2=flatten equation^{3..5};
param=substitute(vars Pb,RP)|substitute(vars Pa,RP);
betti (eq1=substitute(diff(param,transpose eq1),R))
betti (eq1=map(R^{36:2},R^{36:0,78:2},eq1))
eq2=substitute(diff(param,transpose eq2),R);
betti (eq2=map(R^{36:3},R^{36:3,78:2},eq2))
betti (eq=eq1||eq2)
isHomogeneous eq
betti (sol=syz (eq,DegreeLimit=>0)) --51
\end{verbatim}

}

\bigskip

{\bf Acknowledgement}  We would  like to thank Marian Aprodu
and especially
Frank Schreyer and Charles Walter for some useful
conversations.

The present research was carried out in the realm  of the DFG
Schwerpunkt "Globale Methoden in der komplexen Geometrie".
The second author also profited from a travel grant  from a 
DAAD-VIGONI
program.



\end{document}